\documentclass[12pt,draftcls,onecolumn]{IEEEtran}

\ifCLASSINFOpdf
   \usepackage[pdftex]{graphicx}
\else
   \usepackage[dvips]{graphicx}
\fi
\usepackage[cmex10]{amsmath}

\usepackage{algorithmic}

\usepackage{array}

\usepackage{mdwmath}
\usepackage{mdwtab}
\usepackage{amssymb}
\usepackage{graphicx}
\usepackage{subfigure}

\newtheorem{rmk}{Remark}

\def\s {\sigma}
\def\t {\theta}
\def\w {\omega}

\begin{document}
%
\title{Charge-Balanced Minimum-Power Controls for Spiking Neuron Oscillators}
%
%
%

\author{Isuru~Dasanayake and Jr-Shin~Li
\thanks{Isuru Dasanayake and Jr-Shin Li are with the Department
of Electrical and Systems Engineering, Washington University, St. Louis,
MO, 63130, USA. e-mail: dasanayakei@seas.wustl.edu, jsli@seas.wustl.edu}}
\maketitle

\begin{abstract}
In this paper, we study the optimal control of phase models for spiking neuron oscillators. We focus on the design of minimum-power current stimuli that elicit spikes in neurons at desired times. We furthermore take the charge-balanced constraint into account because in practice undesirable side effects may occur due to the accumulation of electric charge resulting from external stimuli. Charge-balanced minimum-power controls are derived for a general phase model using the maximum principle, where the cases with unbounded and bounded control amplitude are examined. The latter is of practical importance since phase models are more accurate for weak forcing. The developed optimal control strategies are then applied to both mathematically ideal and experimentally observed phase models to demonstrate their applicability, including the phase model for the widely studied Hodgkin-Huxley equations.
\end{abstract}


%
\IEEEpeerreviewmaketitle

\section{Introduction}
\label{sec:intro}
The electrical activity of a nervous system and its ability to respond to external electrical signals have been long-standing subjects of active research. The resulting insights have led to the innovation of therapeutic procedures for a wide variety of neurological disorders. Deep brain stimulation is one such method applying electrical pulses to inhibit pathological synchrony among the neurons \cite{Lozano04} and is clinically approved in many countries for the treatment of Parkinson's disease, essential tremor, and Dystonia \cite{Benabid10, Marks05}. The cardiac pace maker is another example in medical practices that employs electric pulses to stimulate nervous tissues in order to regulate a patient's heart rate \cite{Kalafer73, Kalafer76}. In these and many other neurological applications, the use of low power electrical stimuli is desired because, for example, high power stimuli are harmful to biological tissues and the reduction of power consumption in a neurological implant is essential in order to reduce its sizes and lengthen its lifetime. In addition, it is of clinical importance to ensure that any external inputs, e.g., currents, applied to stimulate neurons are charge-balanced. That is, the net amount of the electric charge injected into a neuron over one oscillation cycle should be kept zero, because high levels of the charge accumulation may trigger irreversible electro-chemical reactions, resulting in damage of neural tissues and corrosion of electrodes \cite{Merrill05}.

Many mathematical models have been developed to capture the periodic activities of neuron oscillators \cite{Hodgkin52, Rose89, Morris81, Taylor98} and a well established example is the phase response curve (PRC), which quantifies the asymptotic phase shift of an oscillator due to an infinitesimal perturbation of its state \cite{Brown04}. A phase model accurately approximates the behavior of the corresponding  full state-space system in the neighborhood of its periodic orbit \cite{Izhikevich07}.  Due to their simplicity, phase models are very popular for modeling and analyzing the dynamics of neuron oscillators. For example, the patterns of synchrony resulting from the dynamics of an arbitrary network of oscillators with weak coupling were analyzed using phase models \cite{Ashwin92,Tass89}, and a chain of coupled phase  oscillators has been used to model the lamprey spinal generator for locomotion \cite{Cohen82}. In these studies, the inputs to the oscillatory systems were initially defined, and the dynamical responses of neuron populations were analyzed in detail. Recently, as an alternative objective, control and dynamical systems approaches have been used to manipulate neural activities in a desired way. For instance, minimum-power controls for spiking neurons at specified time instances were derived for some mathematically ideal phase models \cite{Dasanayake11, Moehlis06} and charge-balanced controls were calculated using a numerical shooting method \cite{Nabi09}. Controllability of a network of neurons described by phase models has also been investigated \cite{Li10}.

In this article, we consider a general phase model and derive charge-balanced minimum-power controls for spiking a neuron oscillator at a desired time instance different from its natural spiking time. Both cases of unbounded and bounded control amplitude are examined. The latter is of fundamental and practical importance since there exist physical limitations on medical equipment and safety margins for neural tissues and, more importantly, because phase models are valid under weak forcing. We show that the bounded optimal control has switching characteristics synthesized by the unbounded optimal control and the given control bound. The developed optimal control strategies are then applied to both mathematically ideal, such as sinusoidal, and experimentally observed, such as Hodgkin-Huxley, PRC's to demonstrate their applicability.  In addition, we characterize the range of possible spiking times with respect to the given control amplitude for several phase models. Moreover, we apply the optimal controls derived from the reduced phase model to the corresponding full state-space model to verify the consistency of these models through the reduction and the robustness of our optimal control techniques. Such an important validation is missing in the literature.

This paper is organized as follows. In Section \ref{sec:mim_power_control}, we present the optimal control problem of spiking a general phase oscillator. We find the charge-balanced minimum-power control for a prescribed spiking time with and without a constraint on the control amplitude by using the maximum principle. In Section \ref{sec:example}, we apply the derived optimal control strategies to several commonly-used phase models and present the optimal solutions and numerical simulations. In particular, we calculate optimal controls for experimentally observed PRC's including Morris-Lecar and Hodgkin-Huxley PRCs. These optimal controls produced by the maximum principle are verified by the Legendre pseudospectral computational method \cite{Ross03}.

\section{Optimal charge-balanced controls for spiking neurons}
\label{sec:mim_power_control}
In systems theory, a nonlinear oscillator is described by a set of ordinary differential equations that has a stable periodic orbit. This system of equations can be reduced to a single first order differential equation, which is valid while the state of the full system remains in a neighborhood of its unforced periodic orbit \cite{Brown04}. This reduction allows us to represent the dynamics of a weakly forced oscillator by a single phase variable that defines the evolution of the oscillation. Consider a time-invariant system $\dot{x}=f(x,I)$, where $ x(t)\in\mathbb{R}^n$ is the state and $I(t)\in\mathbb{R}$ is the control, which has an unforced stable attractive periodic orbit $\gamma(t)=\gamma(t+T)$ homeomorphic to a circle, satisfying $\dot{\gamma}=f(\gamma,0)$. We can represent this system in a phase-reduced form as
\begin{equation}
	\label{dynamics}
	\dot{\theta}=f(\theta)+g(\theta)I(t),
\end{equation}  
where $\theta$ is the phase variable, $f$ and $g$ are real-valued functions, and $I(t)\in\mathbb{R}$ is the control \cite{Brown04,Izhikevich07}. One complete oscillation of the system corresponds to $\theta\in[0,2\pi)$. The function $f$ gives system's baseline dynamics and $g$ is known as the phase response curve (PRC), which describes the infinitesimal sensitivity of the phase to an external control input. In the case of neural oscillators, $I$ represents an external current stimulus and $f$ is referred to the instantaneous oscillation frequency in the absence of any external input, i.e., $I=0$. Neuron spiking occurs when the oscillator evolves through one complete cycle. As a convention, the occurrence of spikes takes place at $\theta=2n\pi$, where $n=0,1,2,\ldots$. We consider spiking a neuron at a prescribed time with a minimum-power stimulus and, furthermore, intend to find a charge-balanced one in order to minimize the side-effects cause by the accumulation of electric charge. The design of such charge-balanced minimum-power current stimuli for spiking neurons gives rise to a constrained optimal steering problem for a single-input nonlinear system of the form
\begin{align*}\label{eq:oc1}
    \min_{I(t)}\quad & \int_0^T I(t)^2dt,\\
    {\rm s.t.}\quad  & \dot{\theta}=f(\theta)+g(\theta)I(t),\\
    &\dot{p} =I(t), \\
    &\theta(0)=0,\quad  \theta(T)=2\pi, \tag{P} \\
    &p(0)=0,\quad  p(T)=0,\\
    &|I(t)|\leq M,
\end{align*}
where $M\in \mathbb{R}^+$ defines the bound of the control amplitude, and the time-dependent variable $p(t)=\int_0^t{ I(\s)}d\s$, with boundary conditions $p(0)=p(T)=0$, is introduced to accommodate the charge-balanced constraint. In the following, we first consider the case of an unbounded control, namely, $M=\infty$, and then extend the result to the case when the control is bounded.

\subsection{Charge-Balanced Minimum-Power Control with Unconstrained Amplitude}
Relaxing the amplitude constraint by letting $M=\infty$, we apply the maximum principle to characterize the extremal trajectories. The Hamiltonian of the optimal control problem \eqref{eq:oc1} is given by
\begin{equation}
	\label{Hamiltonian}
	H=\lambda_0I^2+\lambda(f(\theta)+g(\theta)I)+\mu I,
\end{equation}
where $\lambda_0,\lambda$, and $\mu$ are Lagrange multipliers associated with the Lagrangian, system dynamics, and the charge-balanced constraint, respectively. Here we consider normal extremals which are found by taking $\lambda_0\neq0$. Note that more specific abnormal extremals found by letting $\lambda_0=0$ can be analyzed according to the expressions and properties of the functions $f$ and $g$. Our derivations here are made for the general phase model, and therefore these extraordinary cases are omitted. Abnormal extremals are in general uncommon in phase models, and none of the phase models considered in this paper has an abnormal extremal (see Remark \ref{rmk:abnormal} in Section \ref{sec:sinusoidal}). Therefore, without loss of generality, we let $\lambda_0=1$. The optimality condition from the maximum principle demands that $\frac{\partial H}{\partial I}=0$ along the optimal trajectory, which yields  
\begin{equation}
	\label{I}
	I=-\frac{\lambda g(\theta)+\mu}{2}.
\end{equation}
The adjoint variables $\lambda$ and $\mu$ are solutions to the time-varying differential equations $\dot{\lambda}=-\frac{\partial H}{\partial \theta}$ and $\dot{\mu}=-\frac{\partial H}{\partial p}$. Together with \eqref{I} these equations can be written as
\begin{align}
	\dot{\lambda}=&-\lambda \frac{\partial f(\theta)}{\partial\theta}+\frac{\lambda(\lambda g+\mu)}{2}\frac{\partial g(\theta)}{\partial\theta}, \label{adj1}\\
	\dot{\mu}=&0, \label{adj2}
\end{align}
which implies that $\mu$ is a constant. In addition, since the Hamiltonian is not explicitly dependent on time, $H$ is a constant along the optimal trajectory. Hence, we let $H=c$, $\forall\, t\in[0,T]$. This can also be seen from the transversality condition of the maximum principle.

It follows that the optimal multiplier $\lambda$ can be found from \eqref{Hamiltonian} by substituting \eqref{I} for $I$. Then, solving for $\lambda$ yields
\begin{equation}
\label{lambda}
	  \lambda = \frac{-\mu g+2f \pm 2\sqrt{f^2-g\mu f-g^2c}}{g^2}.
\end{equation}
Here we will choose the negative square root because the positive case corresponds to a backward evaluation of the phase, which would invalidate the phase model. The phase velocity equation along the optimal trajectory can then be found by using \eqref{lambda}, \eqref{I}, and \eqref{dynamics}, resulting in
\begin{equation}
	\label{theta}
	\dot{\theta}=\sqrt{f^2-g\mu f-g^2c}.
\end{equation}
In addition, substituting \eqref{lambda} into \eqref{I} gives rise to the optimal control $I^*$ in terms of the two constants $\mu$ and $c$,
\begin{equation}
	\label{I*}
	I^* = \frac{-f+\sqrt{f^2-g\mu f-g^2c}}{g}.
\end{equation}
For a given spiking time $T$, the constants $c$ and $\mu$ can be determined from \eqref{theta} by separation of variables together with the charge-balanced constraint written as $\int_0^{2\pi}{\frac{I^*(\theta)}{\dot{\theta}}d\theta}=0$, which yields 
\begin{eqnarray}
	\label{T}
	T=\int_0^{2\pi}{\frac{1}{\sqrt{f^2-g\mu f-g^2c}}}d\theta,
\end{eqnarray}
and
\begin{eqnarray}
	\label{charge}
	\int_0^{2\pi}{\frac{-f+\sqrt{f^2-g\mu f-g^2c}}{g\sqrt{f^2-g\mu f-g^2c}}}d\theta=0.
\end{eqnarray}
Now the optimal control is completely classified by \eqref{I*}, because the constants $\mu$ and $c$ can be found from \eqref{T} and \eqref{charge} for any specified spiking time $T$.

\begin{rmk}
	In the absence of the charge-balanced constraint, corresponding to $\mu=0$, it is sufficient to characterize the optimal control by \eqref{I*} and \eqref{T}.
\end{rmk}

\subsection{Charge-Balanced Minimum-Power Control with Constrained Amplitude}
\label{sec:charge_balanced_bounded}
In practice, the feasible amplitude of the stimulus is limited, and phase models are valid only for weak forcing. Therefore, spiking neurons with controls of bounded amplitude is of practical importance. In this case where we assume that $|I|\leq M$, $\forall\ t\in[0,T]$, the minimum and maximum possible spiking times for a neuron system can be determined. It is easy to see that for a given bound $M>0$, the minimum spiking time is achieved by
\begin{equation}
	\label{eq:I*Tmin}
    I^*_{T_{min}}=\left\{\begin{array}{ll} M, & g(\theta)\geq 0 \\ -M, & g(\theta)<0, 
    \end{array}\right.
\end{equation}
which keeps the phase velocity at its maximum. The minimum spiking time for a given value of $M$, denoted by $T_{min}^M$, is then given by
\begin{align}
	T_{min}^M=&\int\limits_{\theta\in\mathcal{A}}\frac{1}{f(\theta) +g(\theta)M}d\theta+\int\limits_{\theta\in\mathcal{B}} \frac{1}{f(\theta)-g(\theta)M}d\theta,
\label{TminM}
\end{align}
where the sets $\mathcal{A}$ and $\mathcal{B}$ are defined as
\begin{align*}
	\mathcal{A}=&\left\{\theta|\ g(\theta)\geq 0, 0\leq \theta \leq 2\pi\right\},\\
	\mathcal{B}=&\left\{\theta|\ g(\theta)< 0, 0\leq \theta \leq 2\pi\right\}.
\end{align*} 
Symmetric to the minimum spiking time, the maximum spiking time, denoted $T_{max}^M$, for the bound $M$ can be found by applying the opposite control $-I^*_{T_{min}}$, for $M<\min\{|\frac{f(\theta)}{g(\theta)}|: \, \theta\in[0,2\pi)\}$, and it is given by $T_{max}^M=T^{-M}_{min}$. Note that arbitrarily large spiking times are achievable if the bound $M\geq \min\{|\frac{f(\theta)}{g(\theta)}|: \, \theta\in[0,2\pi)\}$.

It is obvious that if $|I^*(\theta)|\leq M$, $\forall\, \theta\in[0,2\pi)$, then the amplitude constraint is inactive and $I^*$ as in \eqref{I*} is the charge-balanced minimum-power control. While $|I^*|>M$ for some $\theta\in[0,2\pi)$, it is sufficient to consider the case when $I^*>M$ because the case $I^*<-M$ is symmetric. Suppose that $I^*>M$ for $\theta\in(\theta_1,\theta_2)\subset [0,2\pi)$, we now show that the bang control $I=M$ is optimal for $\theta\in[\theta_1,\theta_2]$. Since the Hamiltonian \eqref{Hamiltonian} is a convex function of $I$, $I=M$ is then the minimizer when $I^*>M$ for $\t\in[\t_1,\t_2]$. In this case, we have, from \eqref{Hamiltonian}, the Lagrange multiplier
\begin{equation}
\label{eq:lambdac}
\lambda=\frac{c-M^2-\mu M}{f(\theta)+Mg(\theta)},
\end{equation}
which satisfies the adjoint equation \eqref{adj1}, and hence $I(\theta)=M$ is optimal for $\theta\in[\theta_1,\theta_2]$. Similarly, the same approach can be used to show that $I=-M$ is optimal on the interval over which $I^*<-M$. Therefore, the charge-balanced minimum-power control with limited control amplitude $M$ is of the form with switching characteristic
\begin{equation}
	\label{eq:bounded_control}
    I^*_{M}(\theta)=\left\{\begin{array}{ll}-M, & I^*(\theta)< -M \\ I^*(\theta), & -M\leq I^*(\theta)\leq M\\ M, & I(\theta)^*> M.
    \end{array}\right.
\end{equation}
The switching phases $\theta\in[0,2\pi)$ such that $I^*(\theta)=-M$ or $I^*(\theta)= M$ can be computed (see the examples in Section \ref{sec:example}) and the required parameter values $\mu$ and $c$ can be calculated according to the specified spiking time and the charge-balanced constraint from the equations
\begin{equation}
\label{eq:bounded_control_T}
T=\int_0^{2\pi}\frac{1}{f(\theta)+g(\theta)I^*_M}d\theta
\end{equation}
and
\begin{equation}
\label{eq:bounded_control_CB}
0=\int_0^{2\pi} \frac{I^*_M}{f(\theta)+g(\theta)I^*_M}d\theta.
\end{equation}

\section{Example}
\label{sec:example}
We now apply our optimal control strategies to several commonly-used phase models characterized by various PRC's, including mathematically ideal models, such as sinusoidal PRC, SNIPER PRC, and theta neuron PRC, as well as more realistic phase models such as Hodgkin-Huxley and Morris-Lecar PRC's. These mathematically ideal phase models are approximations to full state-space models at certain bifurcation points, whereas Hodgkin-Huxley and Morris-Lecar phase models are obtained numerically by perturbing their periodic orbits using unit impulses.

\subsection{Sinusoidal Phase Model}
\label{sec:sinusoidal}
The sinusoidal phase model is characterized by a sinusoidal PRC \cite{Brown04},
\begin{eqnarray}
	\label{eq:sin_phasemodel}
	\dot{\t}=\w+z_d (\sin\t) I,
\end{eqnarray}
where $\omega$ is the natural oscillation frequency of the system, $z_d$ is a model-dependent constant, and $I$ is the external stimulus. This is a type II PRC, with both positive and negative regions, which results from a periodic orbit near the supercritical Hopf bifurcation \cite{Brown04}, and occurs in neuron models such as the abstracted FitzHugh-Nagumo neuron model \cite{Keener98}.  Neurons described by this phase model spike periodically with the natural period $T_0=2\pi/\omega$ in the absence of any external input. 

Observe from \eqref{I*} that with $f$ and $g$ as defined above, $I^*(\theta)$ is anti-symmetric around $\theta=\pi$, namely, $I^*(\theta)=-I^*(\theta+\pi)$ for $0\leq\theta\leq\pi$. Therefore, the charge-balanced constraint is automatically fulfilled for the sinusoidal phase model. As a result, we let $\mu=0$.

\subsubsection{Unbounded Control for Sinusoidal Phase Model}
\begin{figure}[ht]
\centering
\begin{tabular}{cc}
    \subfigure[]{\includegraphics[scale=0.18]{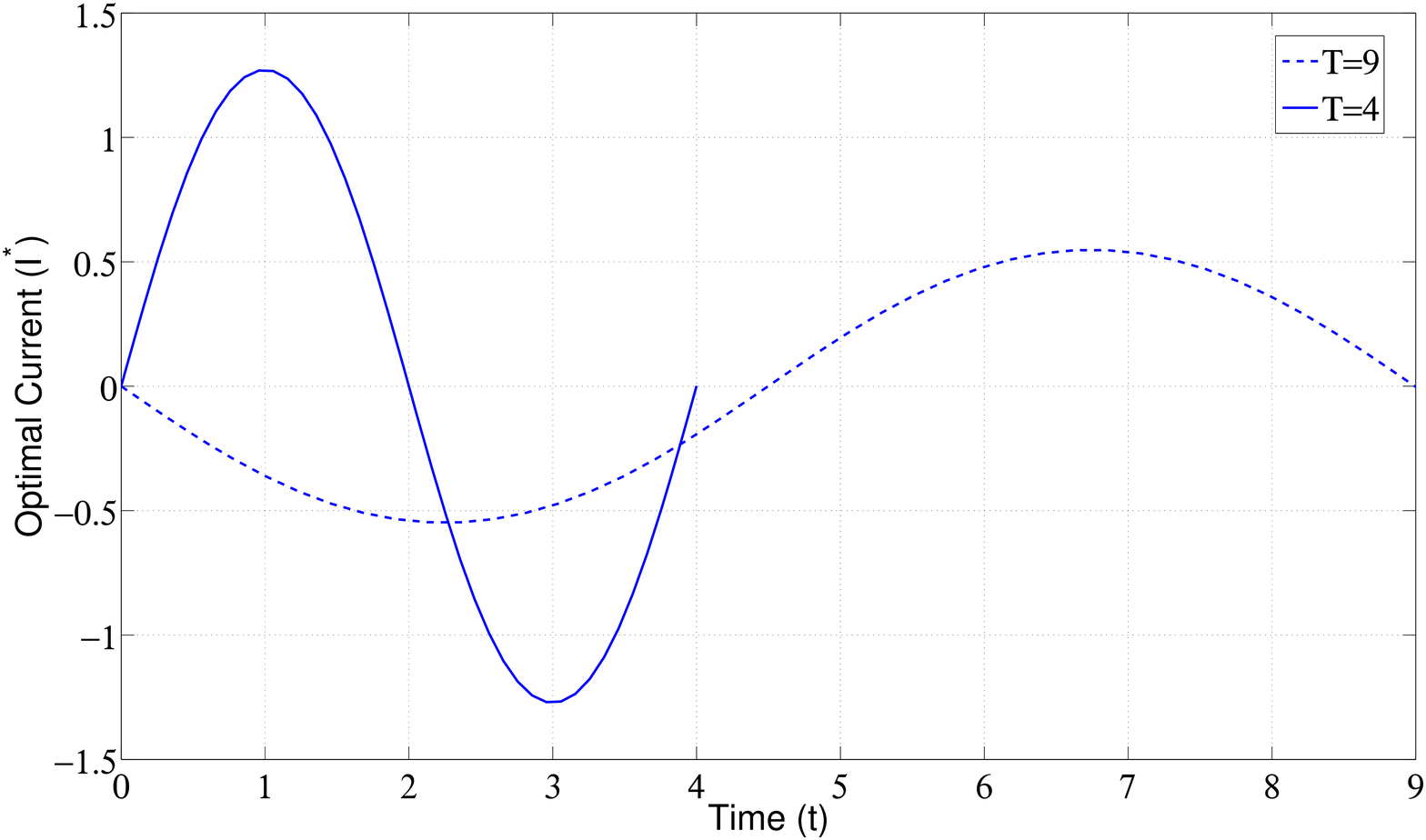} \label{fig:unbounded_control_sine}}
    \subfigure[]{\includegraphics[scale=0.18]{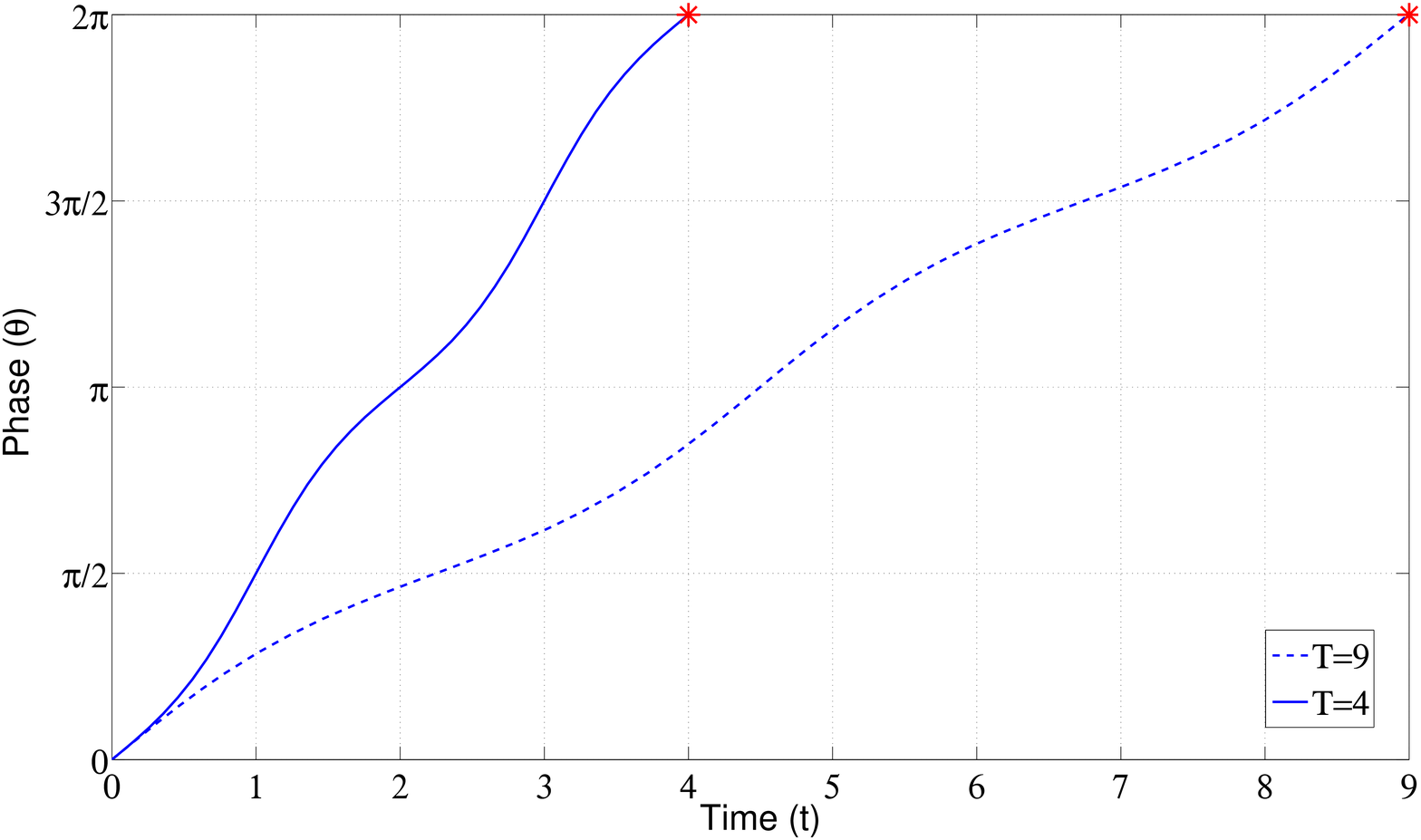} \label{fig:unbounded_trajectories_sine}} &   
\end{tabular}
\caption{Sinusoidal phase model with $\omega=1$ and $z_d=1$. \subref{fig:unbounded_control_sine} Unbounded charge-balanced minimum-power controls for the spiking times $T=4$ and $T=9$. \subref{fig:unbounded_trajectories_sine} Optimal phase trajectories following the optimal controls.}
\end{figure}

Substituting $f=\omega$ and $g=z_d\sin\theta$ with $\mu=0$ into \eqref{I*} and \eqref{T}, the optimal control for spiking a sinusoidal neuron at time $T$ is
\begin{equation}
\label{eq:sine_I*}
	I^* = \frac{-\omega+\sqrt{\omega^2-cz_d^2\sin^2\theta}}{z_d\sin\theta},
\end{equation}
where the constant $c$ is specified by the desired spiking time
\begin{equation}
\label{eq:sine_T}
	T = \int_0^{2\pi}{\frac{1}{\sqrt{\omega^2-cz_d^2\sin^2\theta}}}d\theta.
\end{equation}
A simple example is used to demonstrate these results. For a neuron with the natural oscillation frequency $\omega=1$ and $z_d=1$, the optimal controls for the desired spiking times $T=4$ and $T=9$, smaller and greater, respectively, than the natural spiking time $T_0=2\pi$ are shown in Fig. \ref{fig:unbounded_control_sine}. The corresponding optimal phase trajectories are depicted in Fig. \ref{fig:unbounded_trajectories_sine}.

\begin{rmk}
\label{rmk:abnormal}
	Abnormal extremals in general do not exist in phase models. Consider the case of abnormal extremals for the sinusoidal phase model, where the multiplier $\lambda_0=0$. Then, the Hamiltonian as in \eqref{Hamiltonian} is given by $H=\lambda \omega+ \lambda z_d(\sin\theta)I+ \mu I$, and the optimality condition of the maximum principle gives 
	\begin{equation}
		\label{eq:abnormal_sine}
			\frac{\partial H}{\partial I}=\lambda z_d\sin\theta+\mu=0.
		\end{equation} 
Differentiating this equation with respect to time, we obtain
	\begin{equation}
	\label{eq:abnormal_sine_time_derivative}
	\lambda z_d(\cos\t)\dot{\theta}+\dot{\lambda}z_d\sin\theta+\dot{\mu}=0.
	\end{equation}
	Substituting \eqref{eq:sin_phasemodel}, \eqref{adj1}, and \eqref{adj2} into \eqref{eq:abnormal_sine_time_derivative} for $\dot{\theta},\dot{\lambda}$, and $\dot{\mu}$, respectively, yields
		\begin{equation}
	\label{eq:abnormal_sine_trj}
	\omega\lambda z_d\cos\theta=0.
	\end{equation}
	Abnormal extremals must satisfy \eqref{eq:abnormal_sine_trj}, and it is clear that \eqref{eq:abnormal_sine_trj} holds only when $\lambda\equiv 0$. This leads to $\mu\equiv 0$ from \eqref{eq:abnormal_sine}, which, together with $\lambda\equiv 0$, violate the nontriviality condition of the maximum principle.
\end{rmk}

\subsubsection{Bounded Control for Sinusoidal Phase Model}
\label{case1}
As presented in Section \ref{sec:charge_balanced_bounded}, with the amplitude constraint $|I(t)|\leq M$, $\forall\, t\in[0,T]$, there exists a range of times at which a neuron can be fired. According to \eqref{TminM} for $z_d>0$, the minimum possible spiking time is given by
\begin{equation*}
	T^M_{min}=2\pi\sqrt{\frac{1}{-z_d^2M^2+\omega^2}}-\frac{4\tan^{-1}\left\{z_dM/\sqrt{-z_d^2M^2+\omega^2}\right\}}{\sqrt{-z_d^2M^2+\omega^2}}.
\end{equation*}
Observe from \eqref{eq:sin_phasemodel} that when $M\geq\omega/z_d$, arbitrarily large spiking times can be achieved by making $\dot{\theta}$ arbitrary close to zero. Therefore, the maximum spiking time $T_{max}^M$ is given by $T_{min}^{-M}$ for $M<\omega/z_d$ and the value of $T_{max}^M$ is infinity for $M\geq\omega/z_d$. It follows that the assignment of the spiking time to any $T\in[T^M_{min},T^M_{max}]$ is feasible with the control $I^*_M$ as in \eqref{eq:bounded_control}. Obviously, if the amplitude of the unbounded optimal control satisfies $|I^*|\leq M$ for all $t\in [0,T]$, or equivalently, $\forall\, \theta\in[0,2\pi)$, then the amplitude constraint is inactive and $I^*$ will be the charge-balanced minimum-power control for this bound $M$. There exists a shortest possible spiking time achievable by $I^*$ under the bound $M$, namely (see Appendix \ref{appd:PMP})
\begin{equation}
\label{eq:sin_TI*min}
	T^{I^*}_{min}=\int_0^{2\pi}{\frac{1}{\sqrt{\omega^2+z_dM(z_dM+2\omega)\sin^2\theta}}}d\theta.
\end{equation}
The maximum spiking time achieved by $I^*$ is $T^{I^*}_{max}=T^{I^*}_{min}|_{M=-M}$ for $M<\omega/z_d$ and $T^{I^*}_{max}=\infty$ for $M\geq\omega/z_d$ (see Appendix \ref{appd:PMP}). Note that $T^M_{min}\leq T^{I^*}_{min}\leq T^{I^*}_{max}\leq T^M_{max}$ and a spiking time $T\in(0,T^M_{min})\cup(T^M_{max},\infty)$ cannot be achieved with the bound $M$. In order to properly classify the feasible spiking ranges and associated controls, we consider the two cases, where $M<\omega/z_d$ and $M\geq\omega/z_d$.

\begin{figure}[ht]
	\centering
		\includegraphics[scale=0.7]{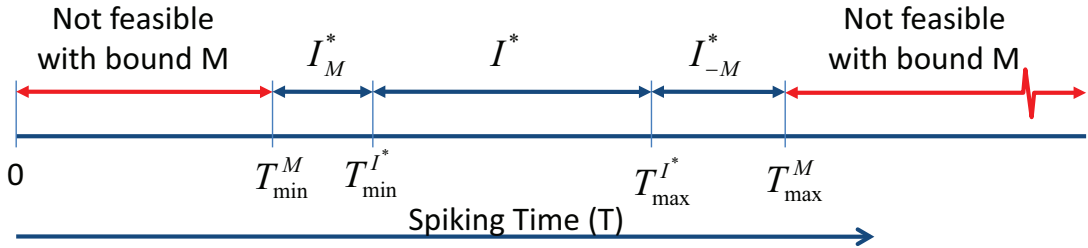}
	\caption{A summary of the optimal control strategies for the sinusoidal PRC model for $M<z_d/\omega$}
	\label{fig:spiking_time_line_case1}
\end{figure}

\begin{figure}[ht]
\centering
    \includegraphics[scale=0.22]{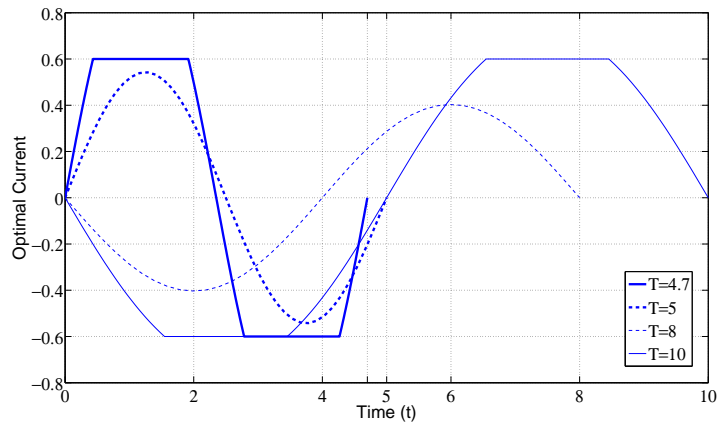} 
\caption{ Optimal bounded controls with bound $M=0.6$ for sinusoidal phase model (with $\omega=1$, $z_d=1$) to elicit spikes at $T=4.7,5,8,10$.}
\label{fig:case2_sine_conl}
\end{figure}

\emph{Case I: ($M<\omega/z_d$})
For a desired spiking time $T\in[T^M_{min},T^{I^*}_{min}]$, the charge-balanced minimum-power control $I^*_M$ is characterized, according to \eqref{eq:bounded_control}, by switching between $I^*$ and $M$,
\begin{equation}
    \label{eq:I1*}
    I^*_{M}=\left\{\begin{array}{ll} I^* & \ 0\leq \theta < \theta_1 \\ M & \ \theta_1\leq\theta\leq\theta_2 \\ I^* & \ \theta_2 < \theta < \theta_3 \\ -M & \ \theta_3\leq\theta\leq\theta_4\\ I^* & \ \theta_4 < \theta \leq 2\pi,
    \end{array}\right.
\end{equation}
in which $\theta_1=\sin^{-1}[-2M\omega/(z_dM^2+z_dc)]$, $\theta_2=\pi-\theta_1$, $\theta_3=\pi+\theta_1$, and $\theta_4=2\pi-\theta_1$ (see Appendix \ref{appd:PMP}). The constant $c$ can be computed according to the desired spiking time $T$, as in \eqref{eq:bounded_control_T}, through the relation
\begin{equation}
	\label{eq:Tbounded_sine}
	T=\int_0^{\theta_1}{\frac{4}{\sqrt{\omega^2-cz_d^2\sin^2\theta}}}d\theta +\int_{\theta_1}^{\frac{\pi}{2}}{\frac{4}{\omega+z_dM\sin\theta}}d\theta.
\end{equation} 
The spiking time $T\in[T^{I^*}_{min},T^{I^*}_{max}]$ can be optimally achieved by the control $I^*$, and for $T\in[T^{I^*}_{max},T^{M}_{max}]$ the optimal control is given by substituting $M=-M$ in the expressions \eqref{eq:I1*} and \eqref{eq:Tbounded_sine}, i.e., $I^*_{-M}$. A summary of the optimal (minimum-power) spiking scenarios for a prescribed spiking time of a neuron governed by the sinusoidal phase model \eqref{eq:sin_phasemodel} is illustrated in Fig. \ref{fig:spiking_time_line_case1}. Fig. \ref{fig:case2_sine_conl} shows the optimal controls  for spiking a sinusoidal neuron with $\w=1$ and $z_d=1$ at $T=4.7,5,8,10$ with the control bound $M=0.6<\w/z_d=1$. These spiking times are chosen to cover all possible spiking scenarios depicted in Fig. \ref{fig:spiking_time_line_case1}. For this particular example, we select the cases of both $T<T_0=2\pi/\w$ and $T>T_0$, where $T=4.7\in[T^M_{\min},T^{I^*}_{\min}]$, $5,8\in[T^{I^*}_{\min},T^{I^*}_{\max}]$ and $10\in[T^{I^*}_{\max},T^M_{\max}]$.

\begin{figure}[ht]
	\centering
		\includegraphics[scale=0.6]{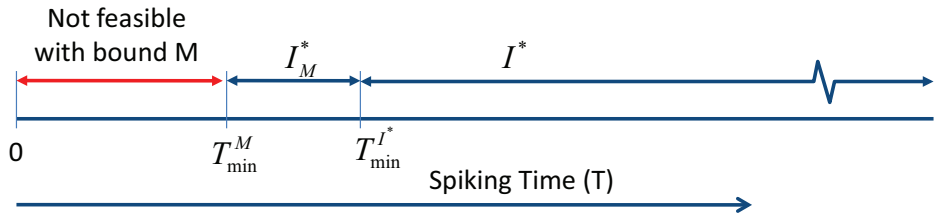}
	\caption{A summary of the optimal control strategies for the sinusoidal PRC model for  $M \geq \omega/z_d$}
	\label{fig:spiking_time_line_case2}
\end{figure}

\begin{figure}[ht]
\centering
 		\includegraphics[scale=0.22]{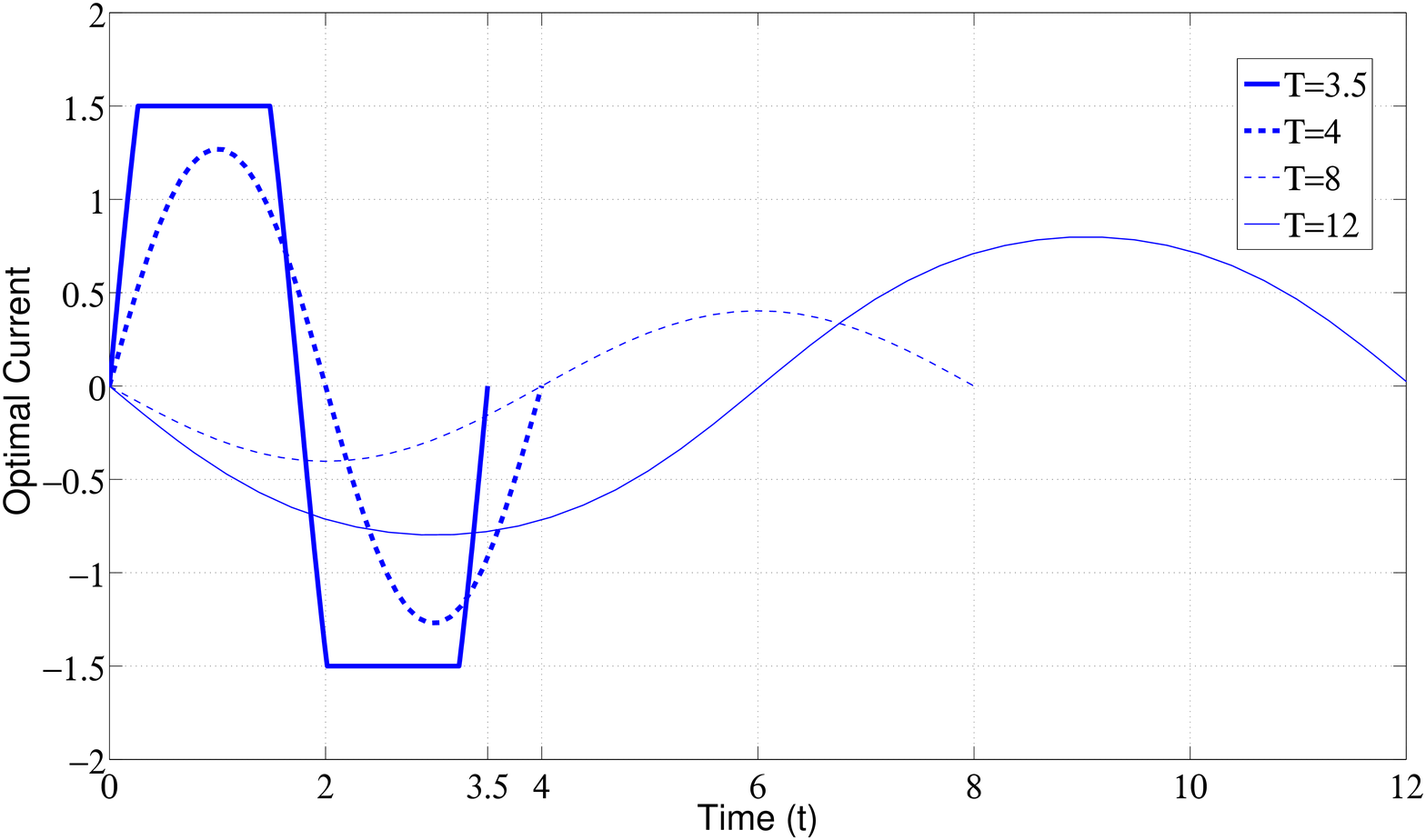} 
\caption{  Optimal bounded controls with bound $M=1.5$ for sinusoidal phase model (with $\omega=1$, $z_d=1$) to elicit spikes at $T=3.5,4,8,12$.}
\label{fig:case1_sine_conl}
\end{figure}

\emph{Case II: ($M\geq\omega/z_d$})
In this case, arbitrarily large spiking times are possible because the system can be driven arbitrarily close to the equilibrium point $\dot{\theta}=0$. Analogous to the previous case, if the desired spiking time is $T\in[T^{M}_{min},T^{I^*}_{min}]$, then the switching control $I^*_M$, as given in \eqref{eq:I1*}, will be optimal, and for $T\in[T^{I^*}_{min},\infty)$ the control $I^*$ will be optimal. A summary of optimal (minimum-power) spiking scenarios for this case is illustrated in  Fig. \ref{fig:spiking_time_line_case2}, and Fig. \ref{fig:case1_sine_conl} shows the optimal controls for spiking a sinusoidal neuron with $\w=1$ and $z_d=1$ at $T=3.5,4,8,12$ given the control bound $M=1.5>\w/z_d=1$. As in the previous case, these spiking times are chosen to cover all possible spiking scenarios depicted in Fig. \ref{fig:spiking_time_line_case2}, for example, $T=3.5\in[T^M_{\min},T^{I^*}_{\min}]$, and $4,8,12\in[T^{I^*}_{\min},\infty]$. Note that in this case $T^{I^*}_{max}=T^{M}_{max}=\infty$.

\subsection{SNIPER Phase Model}
\begin{figure}[ht]
\centering
\begin{tabular}{cc}
		\subfigure[]{\includegraphics[scale=0.18]{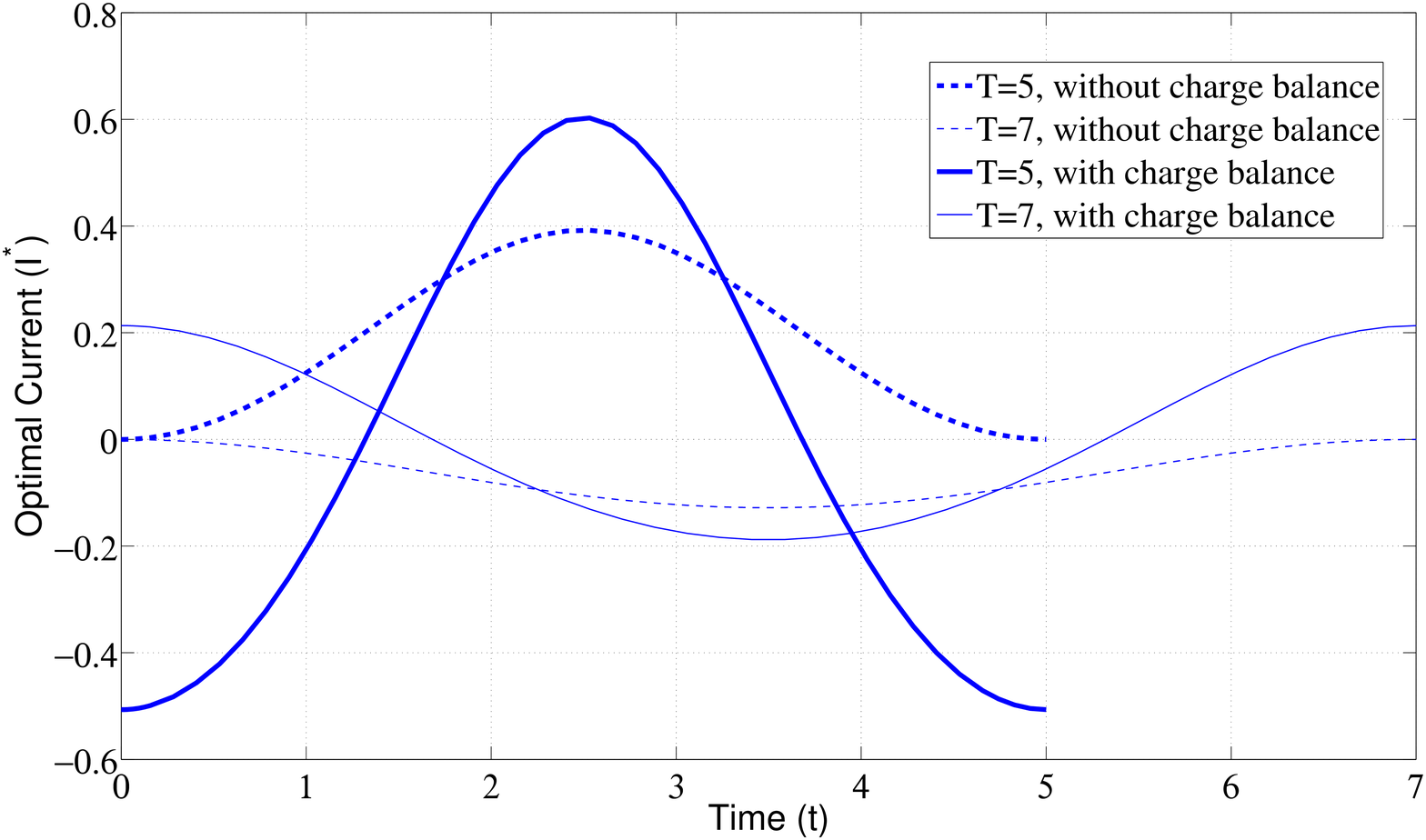} \label{fig:sniper_unbounded_control}} &
    \subfigure[]{\includegraphics[scale=0.18]{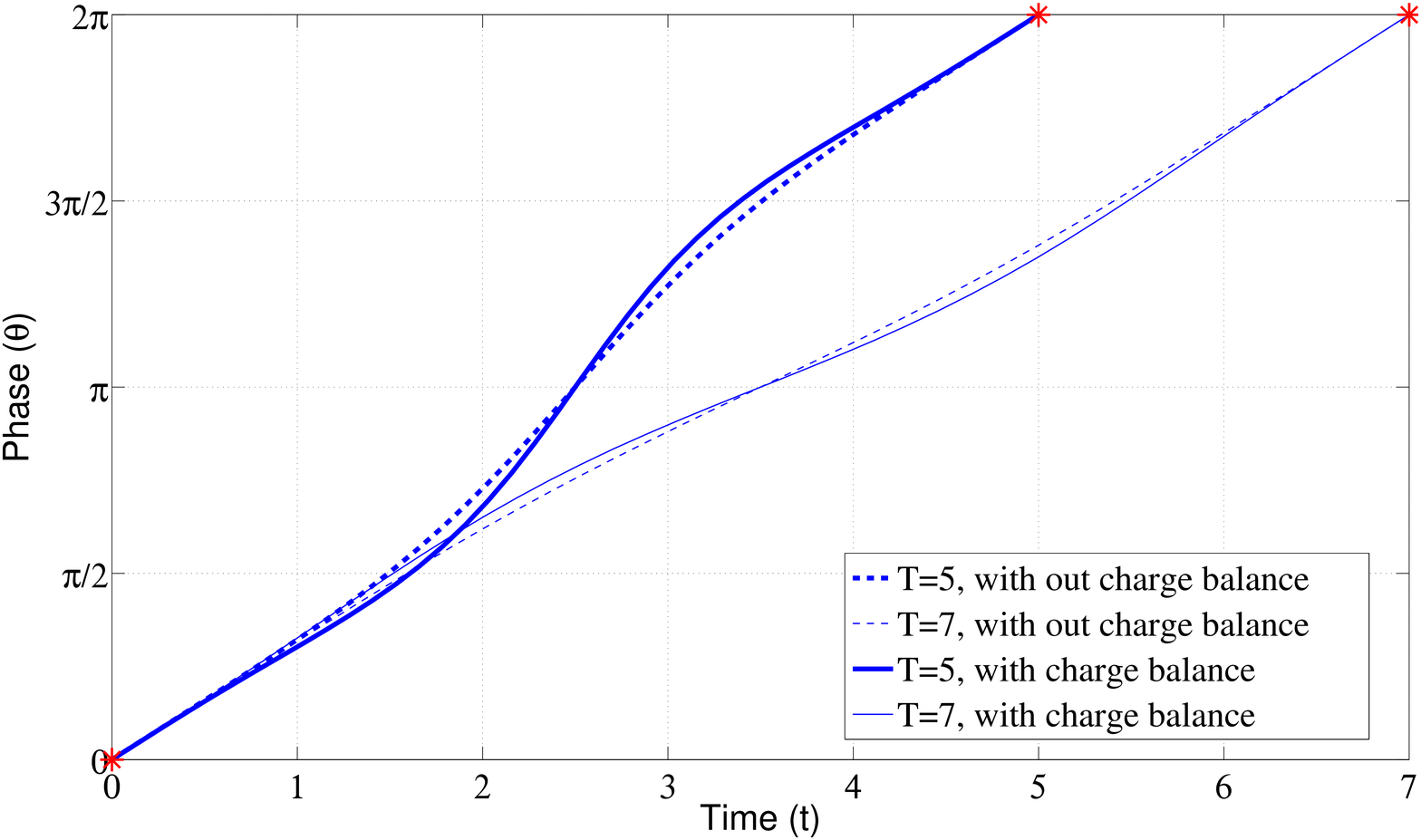} \label{fig:sniper_unbounded_traj}}   
\end{tabular}
\caption{\subref{fig:sniper_unbounded_control} Unbounded optimal controls with and without the charge-balanced constraint for spiking a SNIPER neuron with $\omega=1$ and $z_d=1$ at $T=5$ and $T=7$. \subref{fig:sniper_unbounded_traj} The corresponding optimal phase trajectories.}
\end{figure}

\begin{figure}[ht]
\centering
\includegraphics[scale=0.22]{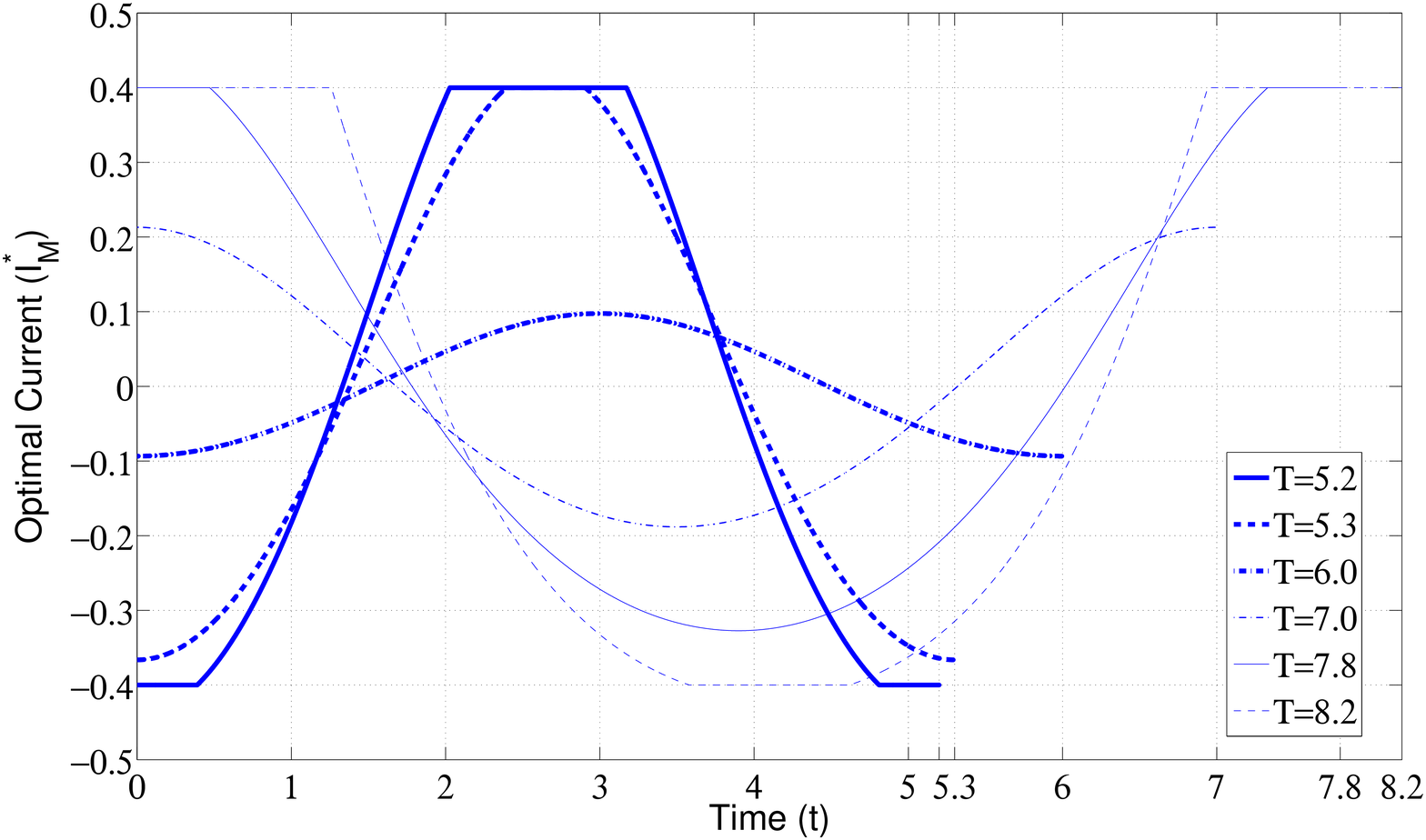} 
\caption{ Optimal charge-balanced controls of minimum power given the control bound $M=0.4$ for spiking a SNIPER neuron with $\omega=1$ and $z_d=1$ at $T=5.2,5.3,6.0,7.0,7.8,8.2$.}
\label{fig:charge_balanced_sniper_bounded_control}
\end{figure}

SNIPER phase model is characterized by $f(\theta)=\omega$ and the PRC $g(\theta)=z_d(1-\cos\theta)$ \cite{Brown04}. This phase model is derived from a SNIPER bifurcation (saddle-node bifurcation of a fixed point on a periodic orbit) which can be found on type I neurons \cite{Ermentrout96} like the Hindmarsh-Rose model. The charge-balanced minimum-power control for unbounded control amplitude can be readily calculated  according to \eqref{I*}, \eqref{T}, and \eqref{charge}, and for bounded control amplitude the control is calculated according to \eqref{eq:bounded_control}, \eqref{eq:bounded_control_T}, and \eqref{eq:bounded_control_CB} using these $f$ and $g$ functions. In Fig. \ref{fig:sniper_unbounded_control} and \ref{fig:sniper_unbounded_traj}, we show unbounded optimal controls in the absence and presence of the charge-balanced constraint and the resulting trajectories of a SNIPER neuron with $\omega=1$ and $z_d=1$. Note that the optimal controls without considering the charge-balanced constraint are obtained by taking $\mu=0$. Fig. \ref{fig:charge_balanced_sniper_bounded_control} illustrates bounded charge-balanced minimum-power controls for spiking the same neuron system at various spiking times which are grater and smaller than its natural spiking period $T_0=2\pi$. We present controls driving the neuron from $\t=0$ to $\t=2\pi$ at various times, $T=5.2$, $5.3$, $6.0$, $7.0$, $7.8$, $8.2$. There exist three structurally different controls which have four switches, two switches, and zero switches, depending on the desired spiking time. 


\subsection{Theta Neuron Phase Model}
\begin{figure}[ht]
\centering
\includegraphics[scale=0.22]{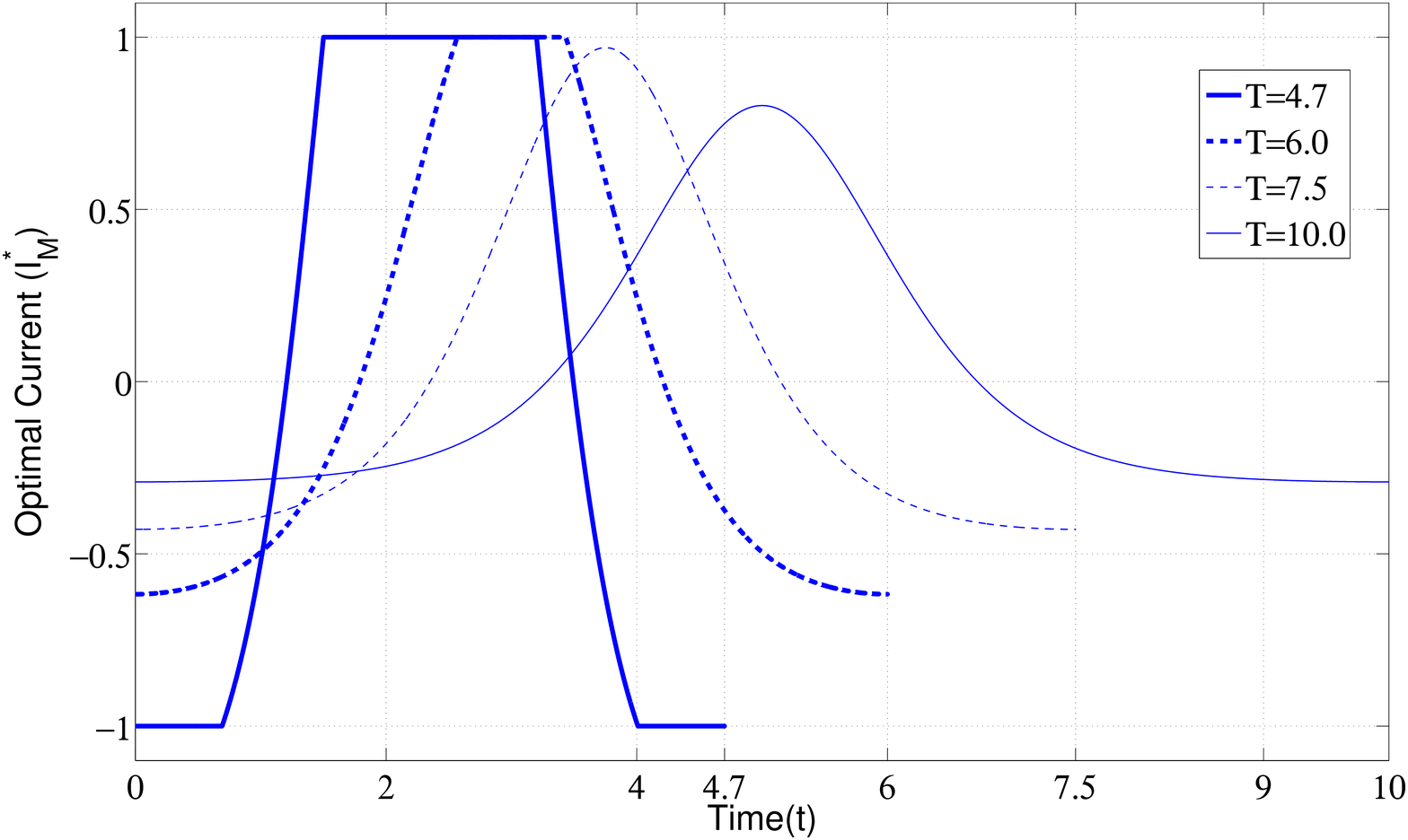}  

\caption{Optimal charge-balanced controls with bound $M=1.0$ and for Theta neuron model (with $I_b=-0.25$) to elicit spikes at $T=4.7,6.0,7.5,10.0$.}
\label{fig:charge_balanced_theta_bounded_control}
\end{figure}

The theta neuron phase model is defined by $f(\theta)=1+\cos\theta+(1-\cos\theta)I_b$ and $g(\theta)=(1-\cos\theta)$, where $I_b$ is known as the neuron baseline current \cite{Moehlis06}. If $I_b>0$, then the neuron spikes with the period $T_0=\pi/\sqrt{I_b}$ in the absence of any external current $I(t)$. When $I_b\leq0$, the neuron does not spike autonomously but it can be fired by the use of an input $I(t)$. Since for $I_b>0$ this neuron model can be transformed to the SNIPER phase model by a coordinate transformation \cite{Moehlis06}, we focus here on the case of $I_b<0$. Similarly, the unbounded and bounded charge-balanced minimum-power controls can be directly calculated by employing \eqref{I*}, \eqref{T}, and \eqref{charge}, or \eqref{eq:bounded_control}, \eqref{eq:bounded_control_T}, and \eqref{eq:bounded_control_CB} in Section \ref{sec:mim_power_control}, respectively. Optimal controls for spiking a theta neuron with $I_b=-0.25$ and $M=1$ are shown in Fig. \ref{fig:charge_balanced_theta_bounded_control}.

The above phase models, though commonly used, are ideal mathematical models of neuron oscillators. We now apply our optimal control strategies to models with experimentally observed PRC's, such as Hodgkin-Huxley and Morris-Lecar phase models, to demonstrate their applicability and generality.

\subsection{Morris-Lecar Phase Model}
\label{sec:ML}
\begin{figure}[ht]
\centering
\begin{tabular}{cc}
    \subfigure[]{\includegraphics[scale=0.18]{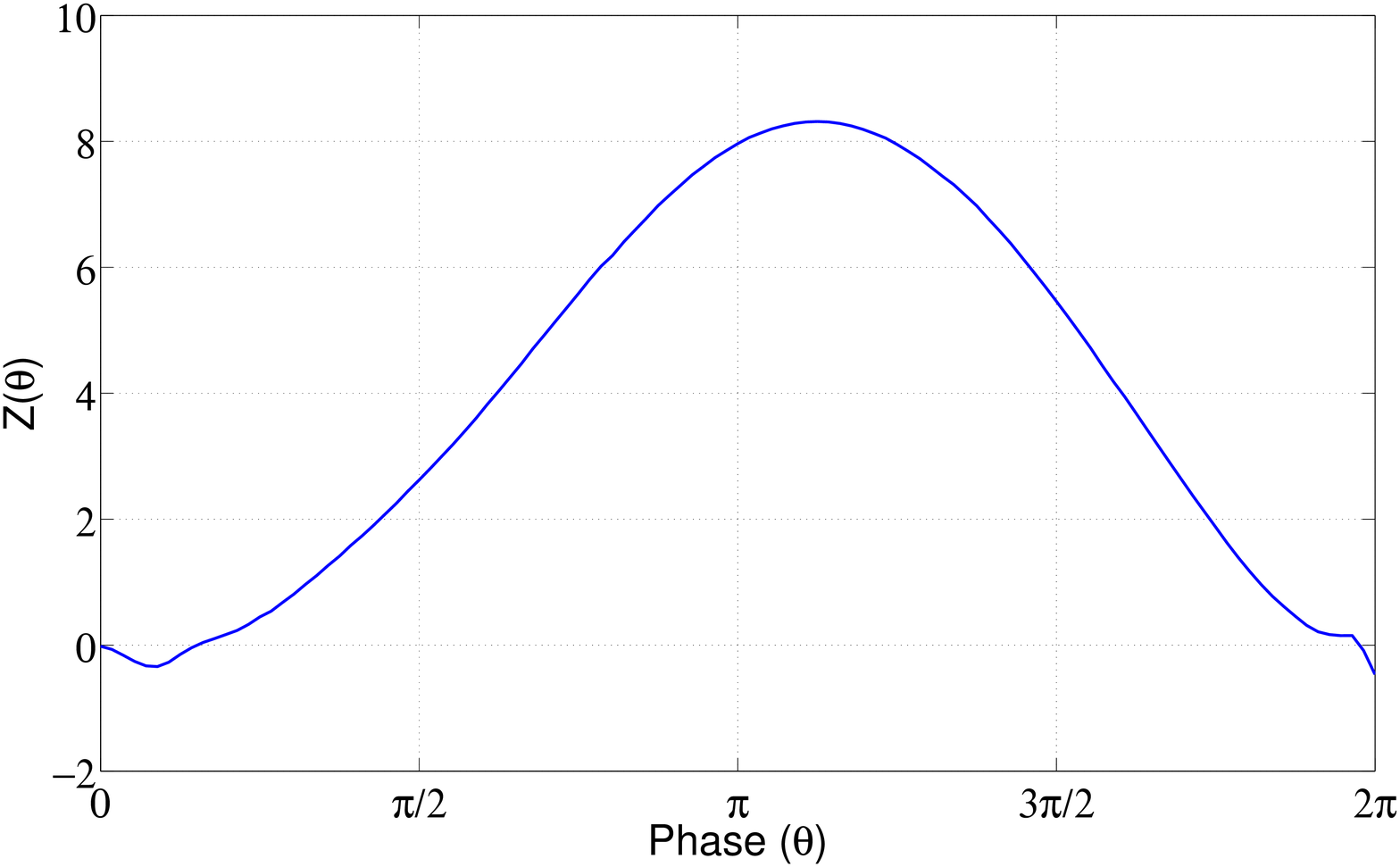} \label{fig:ML_PRC}} &
    \subfigure[]{\includegraphics[scale=0.18]{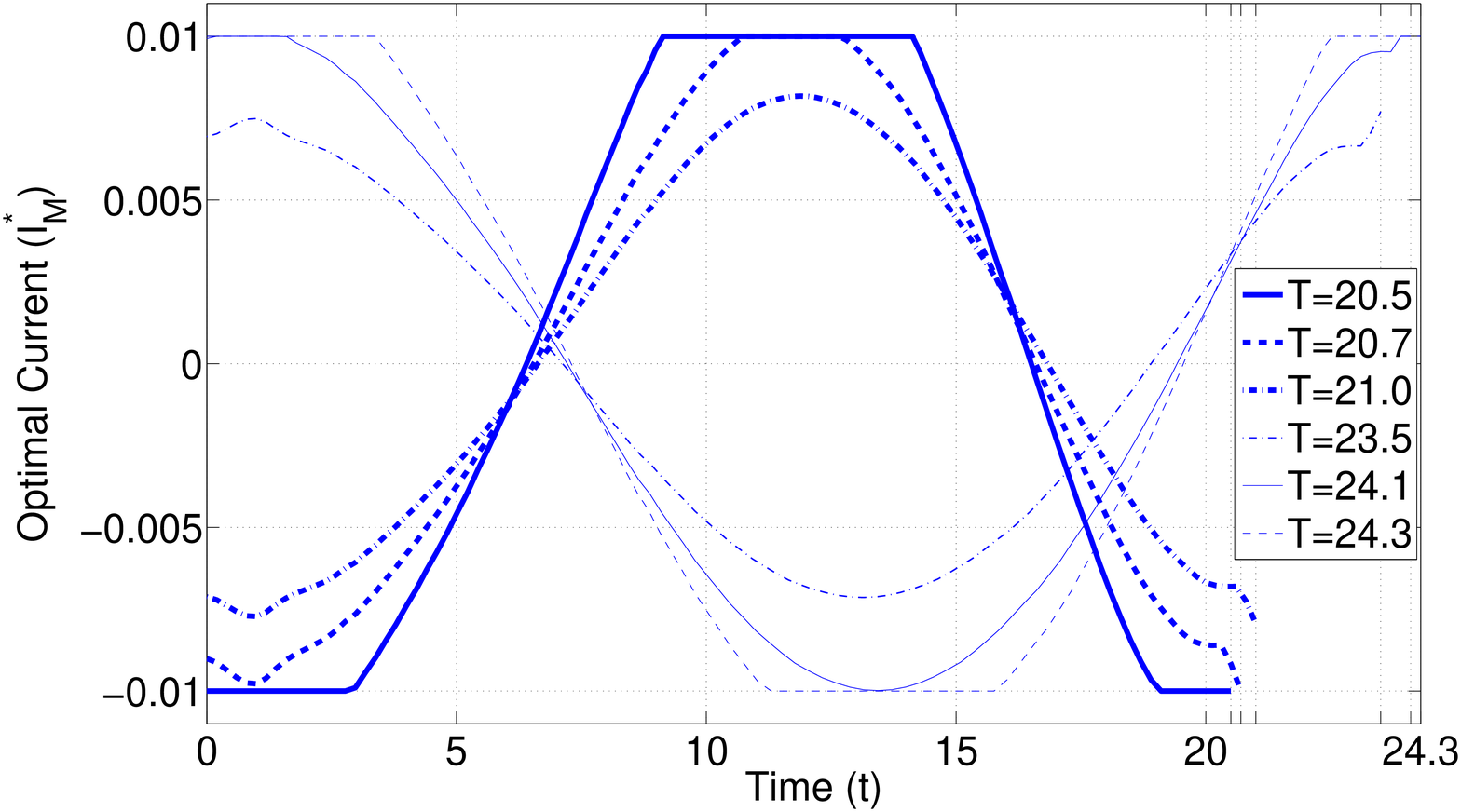} \label{fig:CB_ML_bounded_control}}
\end{tabular}
\caption{\subref{fig:ML_PRC} The Morris Lecar PRC for the parameters given in Appendix \ref{appd:ML}. \subref{fig:CB_ML_bounded_control} Optimal charge-balanced controls of minimum power for spiking a Morris-Lecar neuron at $T=20.5,20.7,21.0,23.5,24.1,24.3 \ ms$ given the control bound $M=0.01\ \mu A$.}
\label{fig:ML_analytical}
\end{figure}

The Morris-Lecar model was originally proposed to capture the oscillating voltage behavior of giant barnacle muscle fibers (see Appendix \ref{appd:ML}) \cite{Morris81}. Over the past years this model has been extensively studied and used as a standard model for representing many different real neurons that are experimentally observable. For example, it has been found that Morris-Lecar PRC is extremely similar to the experimentally observed PRC's of Aplysia motoneuron \cite{Milton00}. The phase model of the Morris-Lecar neuron is given by
\begin{equation}
\label{eq:ML_phase_model}
\dot{\theta}=\omega+Z(\theta)I(t),
\end{equation}
where $\omega$ is the natural oscillation frequency and $Z(\theta)$ represents the PRC which can be calculated numerically from the ODE system in Appendix \ref{appd:ML} by the software package XPP \cite{Ermentrout02}. For the set of parameter values given in Appendix \ref{appd:ML}, the natural frequency $\omega_{ML}=0.283\ rad/ms$ and the PRC is depicted in Fig. \ref{fig:ML_PRC}. The charge-balanced minimum-power controls that elicit spikes for this phase model at various times are shown in Fig. \ref{fig:CB_ML_bounded_control}. We consider six different cases for which the optimal controls have zero, two, and four switchings for spiking times that are longer and shorter than the natural spiking time, $T_0=2\pi/\w_{ML}=22.202 \ ms$.

\subsection{Hodgkin-Huxley Phase Model}
\label{sec:HH}
\begin{figure}[ht]
\centering
\begin{tabular}{cc}
    \subfigure[]{\includegraphics[scale=0.18]{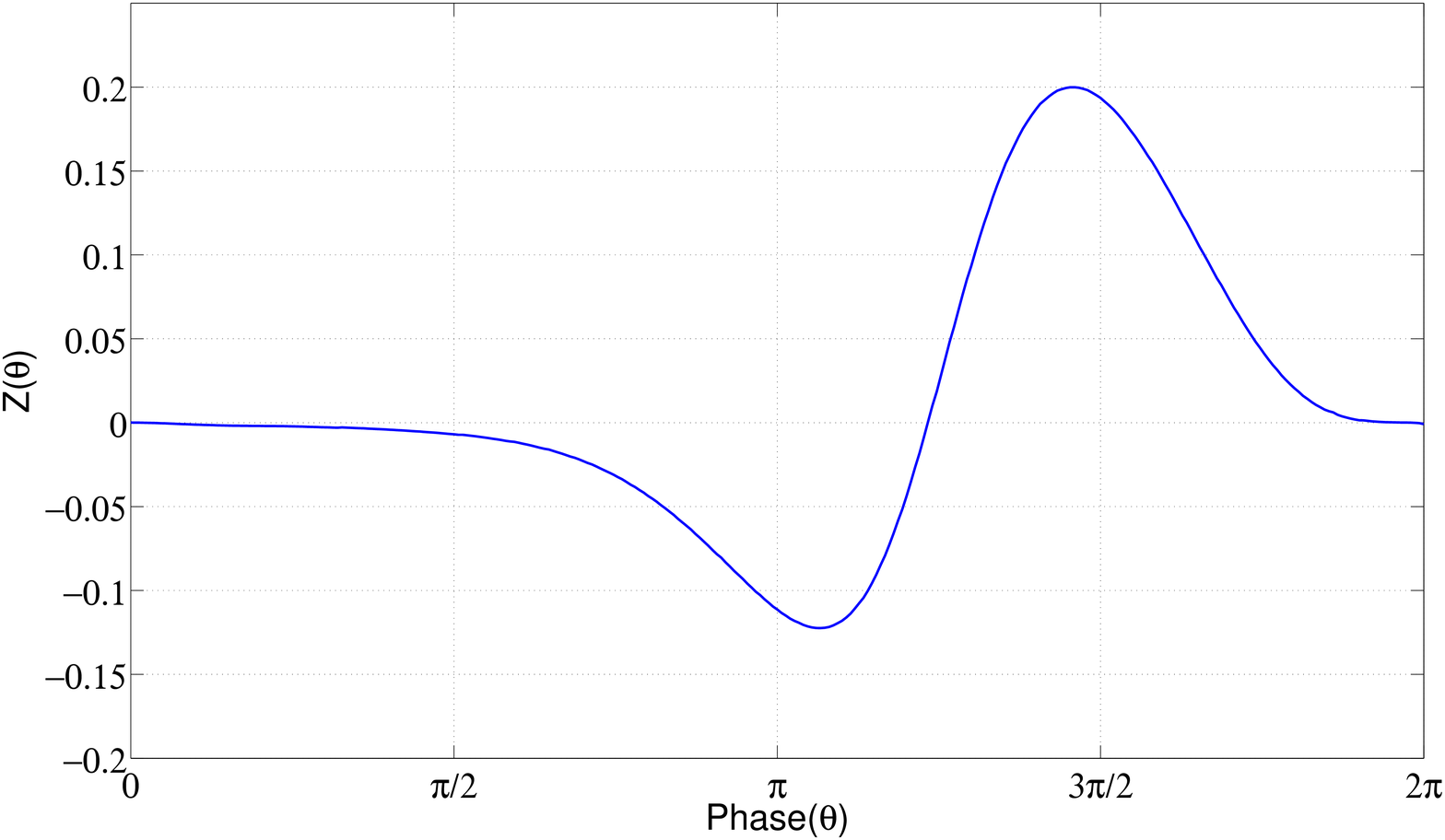} 
	\label{fig:HH_PRC}} &
    \subfigure[]{\includegraphics[scale=0.18]{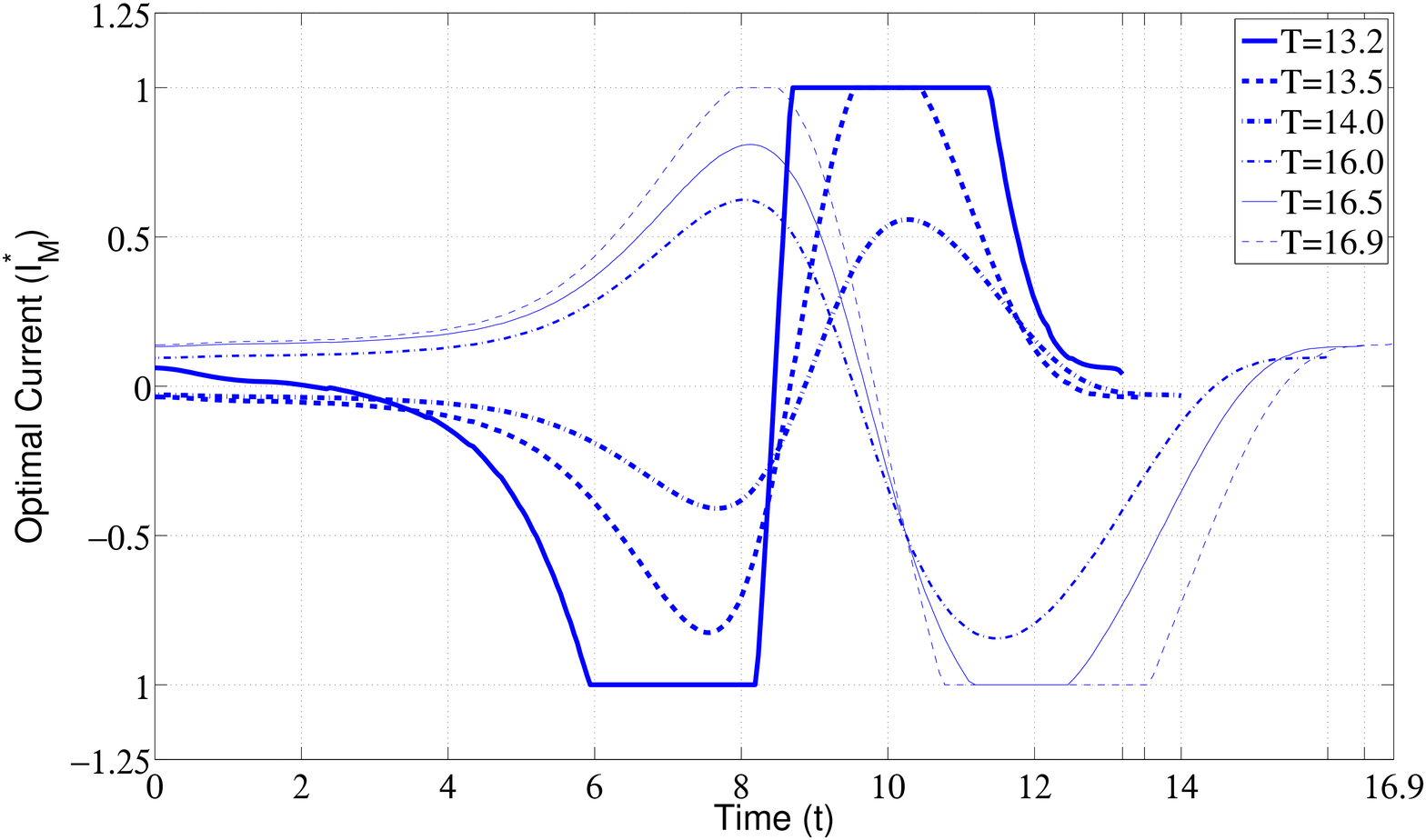}
 	\label{fig:charge_balanced_HH_bounded_control}}
\end{tabular}
\caption{\subref{fig:HH_PRC} The Hodgkin Huxley PRC for the parameters given in Appendix \ref{appd:HH}. \subref{fig:charge_balanced_HH_bounded_control} Optimal charge-balanced controls of minimum power for spiking a Hodgkin-Huxley neuron at $T=13.2,13.5,14.0,16.0,16.5,16.9 \ ms$ given the control bound $M=1.0 \ mA$.}
\label{fig:HH_analytical}
\end{figure}

The Hodgkin-Huxley neuron model is a four dimensional system that describes the propagation and initiation of the action potential in squid axon (see Appendix \ref{appd:HH}) \cite{Hodgkin52}. The phase model for this neuron oscillator is also of the form as in \eqref{eq:ML_phase_model}. For the set of parameter values given in Appendix \ref{appd:HH}, the system has a natural frequency $\w_{HH}=0.4292\ rad/ms$ and its PRC is displayed in Fig. \ref{fig:HH_PRC}. The charge-balanced minimum-power controls that elicit spikes at different time instances are shown in Fig. \ref{fig:charge_balanced_HH_bounded_control}. 

\begin{figure}[ht]
\centering
\begin{tabular}{cc}
		\subfigure[]{\includegraphics[scale=0.18]{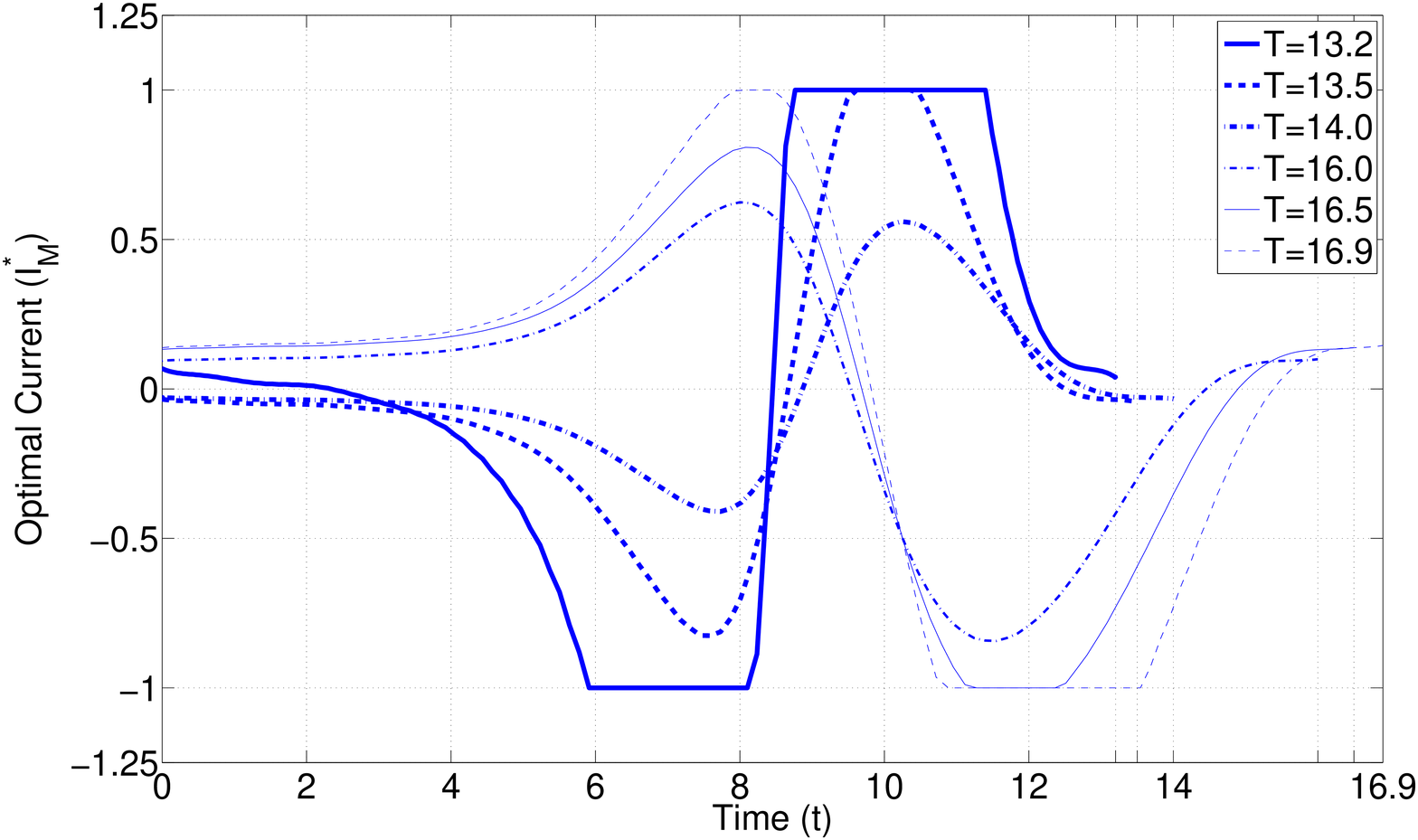} \label{fig:CB_HH_bounded_numerical_controls}} &
    \subfigure[]{\includegraphics[scale=0.18]{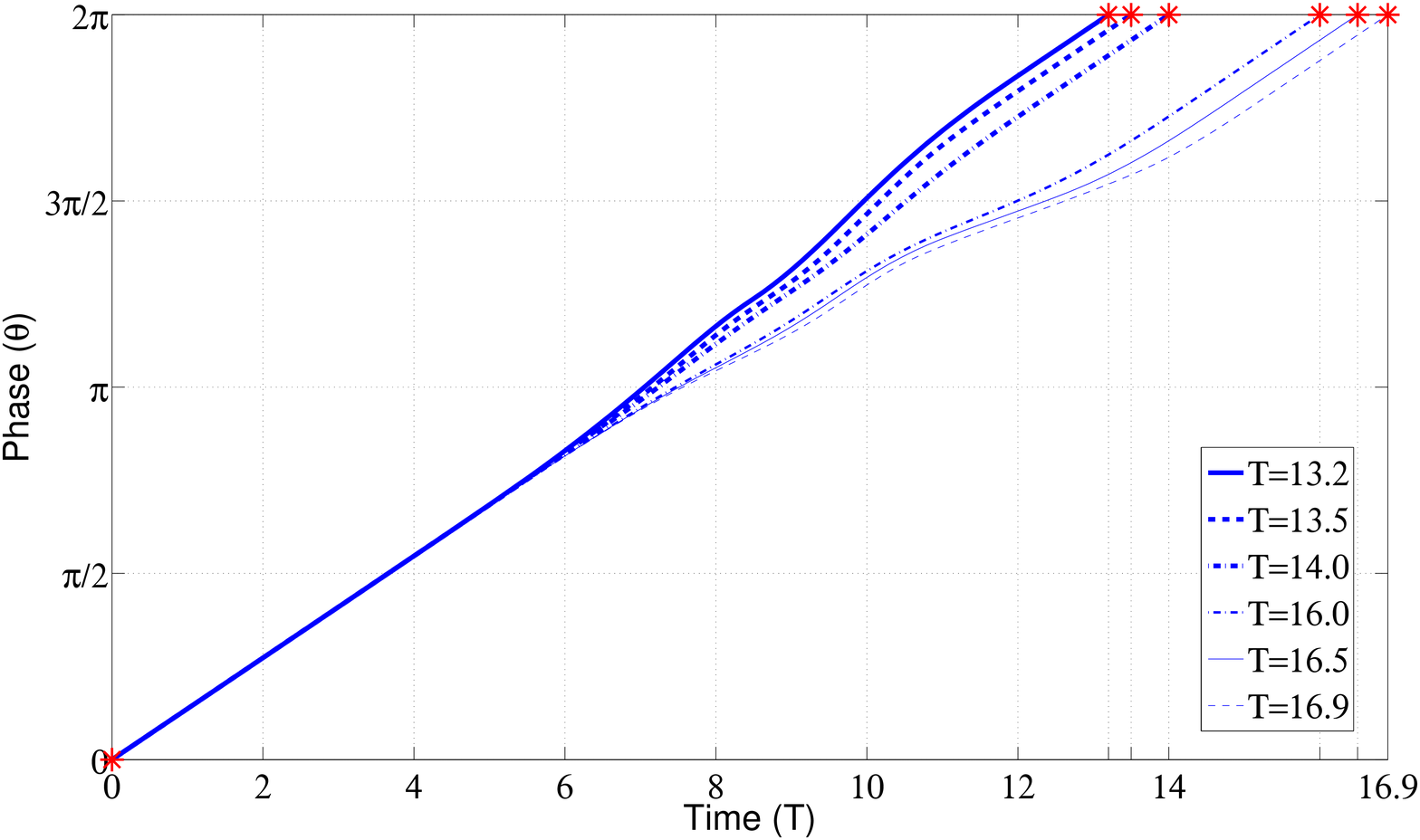} \label{fig:CB_HH_bounded_numerical_traj}} 
\end{tabular}
\caption{\subref{fig:CB_HH_bounded_numerical_controls} Charge-balanced minimum-power controls generated by the Legendre pseudospectral method, which show excellent agreement with the theoretically calculated optimal controls as shown in Fig. \ref{fig:charge_balanced_HH_bounded_control}. \subref{fig:CB_HH_bounded_numerical_traj} The corresponding optimal phase trajectories for the Hodgkin-Huxley neuron.}
\label{fig:HH_numerical}
\end{figure}

Finally, we verified these optimal controls derived with the maximum principle by using the Legendre pseudospectral method. This computational method is a direct and powerful method for solving continuous-time optimal control problems. The basic principle is described in Appendix \ref{sec:PSmethod} and those readers interested in this method can refer to the recent comprehensive work in this area \cite{Gong06, Li09, Li_JCP11, Li_PNAS11}. The optimal controls generated by this pseudospectral method are presented in Fig.\ref{fig:HH_numerical}, which show excellent agreement with the theoretically calculated ones given in Fig.\ref{fig:charge_balanced_HH_bounded_control}.

\begin{figure}
	\centering
		\includegraphics[scale=0.3]{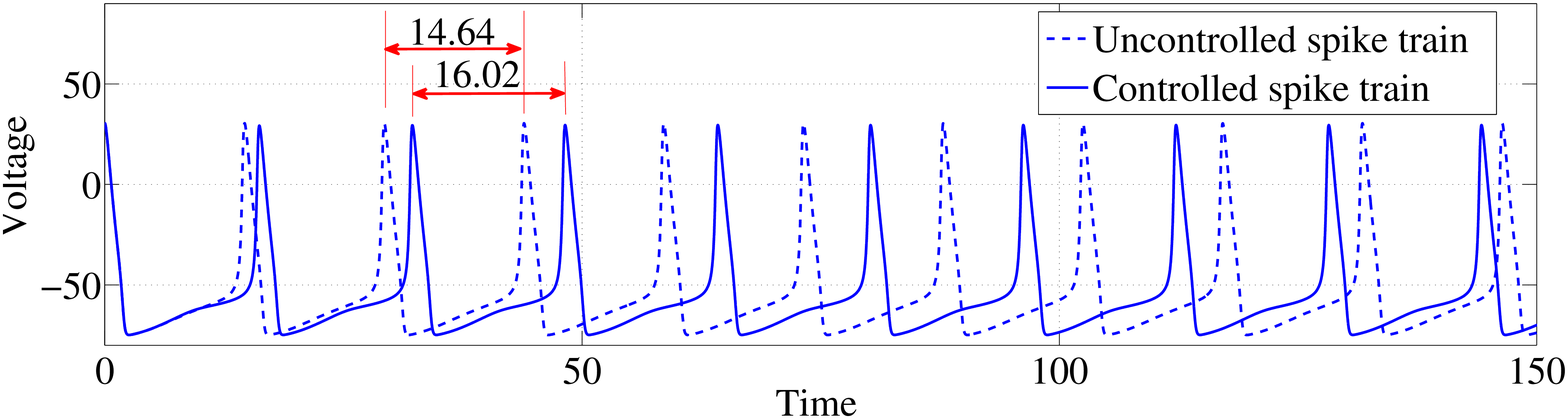}
	\caption{Uncontrolled and controlled spiking trains of Hodgkin-Huxley Model}
	\label{fig:uncontrolled_and_controlled_spike_train_HH}
\end{figure}

Phase models characterize the reduced dynamic behavior of the underlying oscillating systems, where the phase, but not the full state, can be observed. There is a fundamental need to explore the limits of the phase-reduced model as an approximation to the original oscillating system, because this important validation is largely lacking in the literature. The optimal controls for phase models presented so far in this work change the spiking times of an oscillator during the course of one oscillatory cycle, so that a desired spike train can be constructed by repeating the control input. We now apply the optimal controls derived according to the scalar Hodgkin-Huxley phase model to its full state-space model, which is a system of four differential equations as shown in Appendix \ref{appd:HH}. The spike train obtained by repeated application of the optimal control producing an inter-spike time $T=16\ ms$, subject to the control amplitude bound $M=1\ mA$, and the uncontrolled train spiking at the natural period, $T_0=2\pi/\w_{HH}=14.64\ ms$, are illustrated in Fig. \ref{fig:uncontrolled_and_controlled_spike_train_HH}. It is seen that the optimal control delays the spiking time from $14.64\ ms$ to $16.02\ ms$
in the state-space model.

\section{Discussion}
\label{sec:discussion}
In this paper, we considered the optimal control of phase models of neuron oscillators. We derived charge-balanced minimum-power current stimuli that elicit spikes of neurons at desired time instances for the cases of unbounded and bounded control amplitude. In particular, we showed that for the bounded case the optimal control has switching characteristics synthesized by the unbounded optimal control and the control bound. We implemented the resulting analytical optimal controls to various commonly used phase models, including mathematically ideal and experimentally observed models, to demonstrate their applicability. We then applied the optimal controls derived according to the phase-reduced model of Morris-Lecar and Hodgkin-Huxley to the corresponding full state-space system to validate the approximation of the phase model under weak forcing. The theory presented in this work can be applied not only to neuron oscillators but also to any oscillating systems that can be represented using similar model reduction techniques such as biological, chemical, electrical, and mechanical oscillators.


The theoretical results presented in this paper characterize the fundamental limit of how the dynamics of neurons can be perturbed by the use of external inputs. Alternatively, they provide an insight into how the neuron dynamics determine the synaptic input necessary for eliciting spikes, which facilitates the development of optimal stimuli for neurological treatments such as deep brain stimulation for Parkinson's disease. The extension of this work to the optimal control of networks of neuron oscillators is of fundamental and practical importance. Our recent work has shown that a ensemble of uncoupled neurons is controllable, and that the minimum-power controls that spike a network of heterogeneous neurons can be found by using a multidimensional pseudospectral method \cite{Li10}. We plan to extend this recent work to investigate controllability of coupled neurons and related optimal control problems. Systems described by the Kuramoto model will be considered.


%



\appendices
\section{Spiking Sinusoidal Neurons with Bounded Control}
\label{appd:PMP}
\begin{figure}[ht]
\centering
\begin{tabular}{cc}
	\subfigure[]{\includegraphics[scale=0.18]{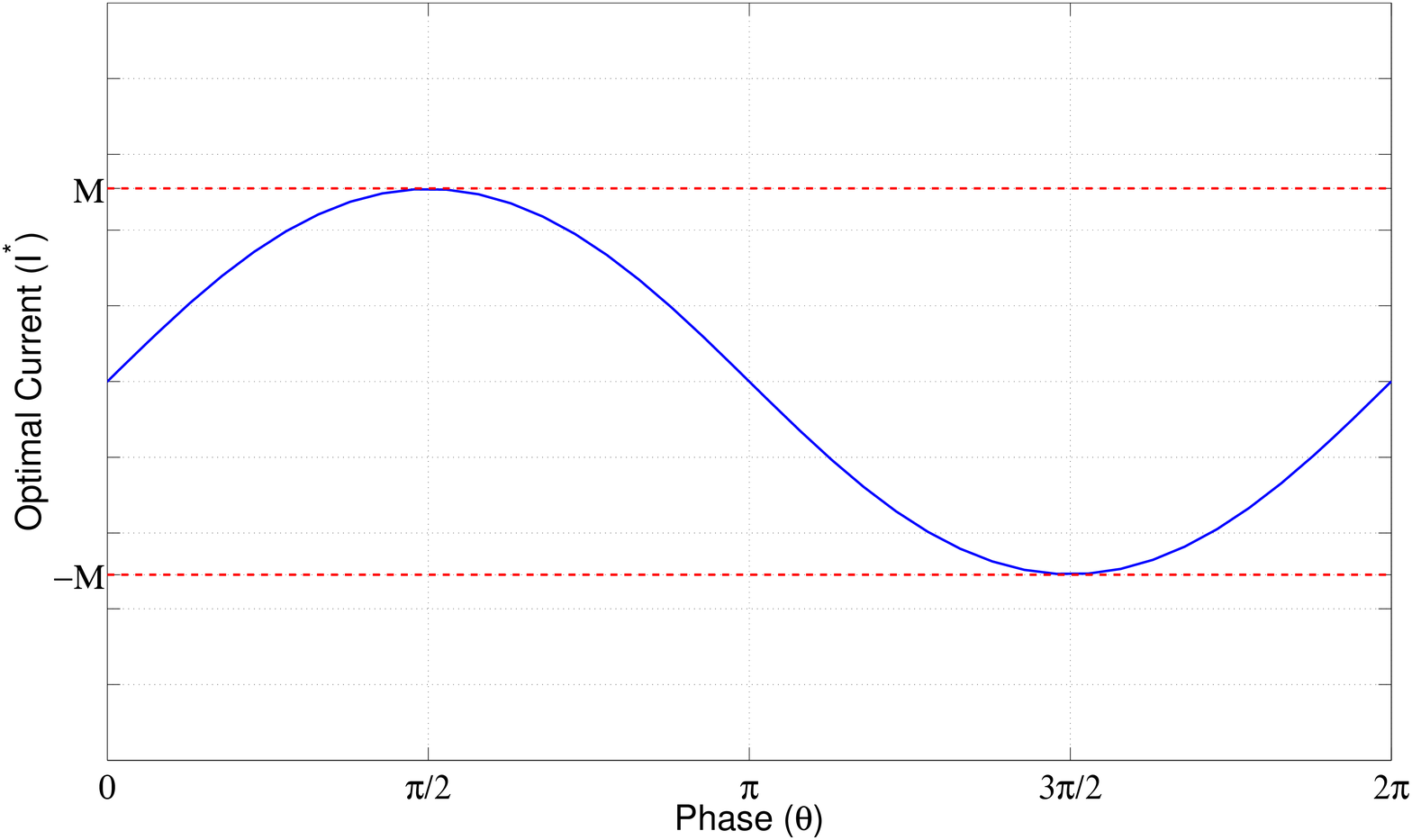} \label{fig:TImin_cal_sine}} &
    \subfigure[]{\includegraphics[scale=0.18]{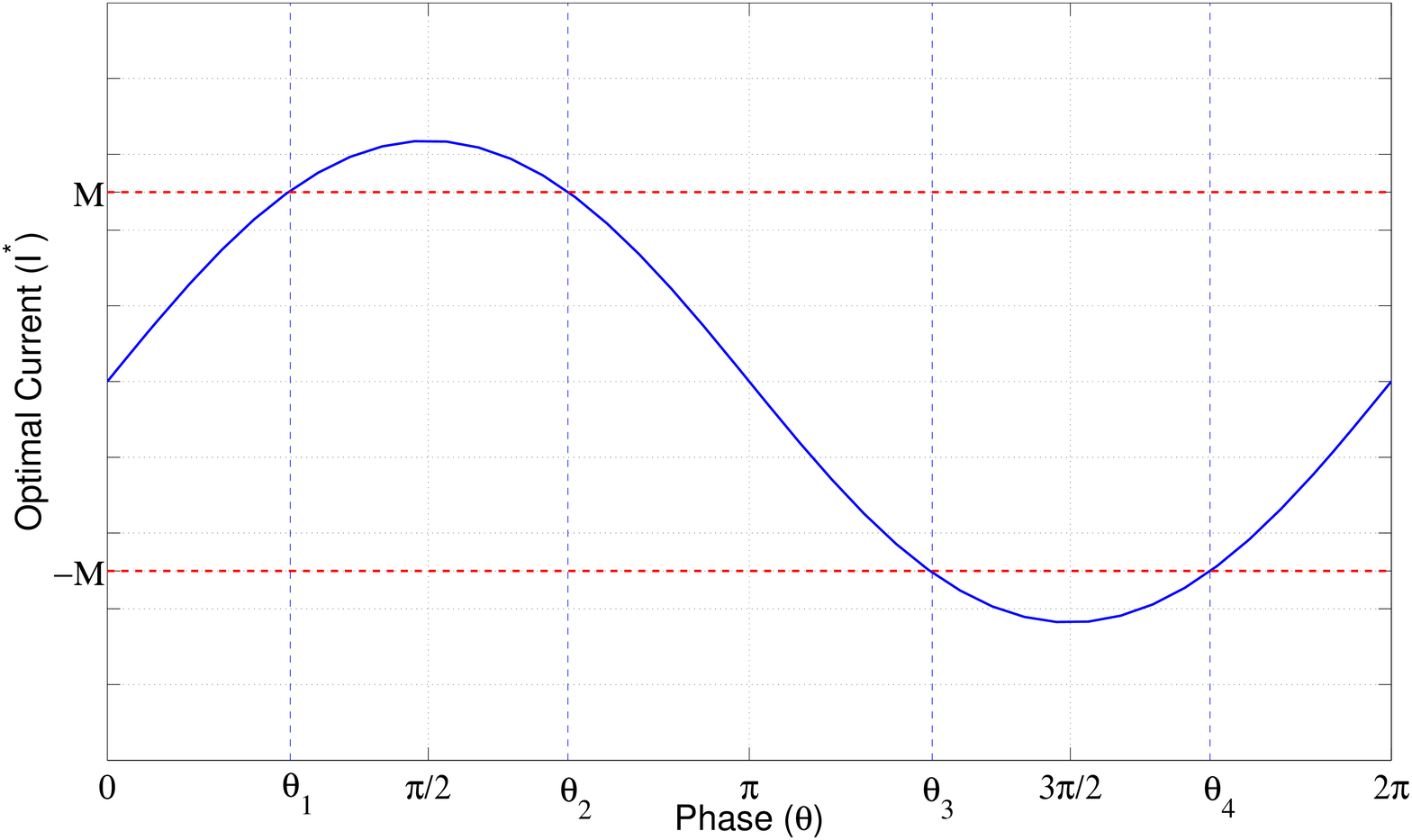} \label{fig:sw_cal_sine}}   
\end{tabular}
\caption{\subref{fig:TImin_cal_sine} An illustration of the optimal control $I^*$ with its maximum value occurring at $\t=\pi/2$ for $c>0$, which gives the shortest possible spiking time subject to the control bound $M$. \subref{fig:sw_cal_sine} An illustration of the case when $I^*>M$ with intersections at $\t_1$, $\t_2$, $\t_3$, and $\t_4$.}
\end{figure}

Simple first and second order optimality conditions applied to \eqref{eq:sine_I*} find that the maximum value of $I^*$ occurs at $\theta=\pi/2$ for $c<0$ and at $\theta=3\pi/2$ for $c>0$ (see Fig.\ref{fig:TImin_cal_sine} for $c<0$). According to \eqref{eq:sine_T}, $c=0$ corresponds to $T=2\pi/\omega$ and $c<0$ ($c>0$) corresponds to $T<2\pi/\omega$ ($T>2\pi/\omega$). Therefore, the constant $c$ for the shortest spiking time with the control $I^*$ satisfying $|I^*(t)|\leq M$ can be calculated by substituting $I^*=M$ and $\theta=\pi/2$ to \eqref{eq:sine_I*}, and then from \eqref{eq:sine_T} we obtain the shortest spiking period by $I^*$, $T_{min}^{I^*}$, as in \eqref{eq:sin_TI*min}.

Since $I^*$ takes the maximum value at $\t=3\pi/2$ for $c>0$, which corresponds to $T>2\pi/\omega$, we have $|I^*|\leq(\omega-\sqrt{\omega^2-cz_d^2})/z_d$, which leads to $|I^*|<\w/z_d\leq M$ for $T>2\pi/\omega$, provided that $M\geq\w/z_d$. This implies that $I^*$ is the minimum-power control for any desired spiking time $T>2\pi/\w$ when $M\geq\w/z_d$. Since the smallest spiking time by the control $I^*$ with the bound $M$ is given by $T^{I^*}_{min}$ as in \eqref{eq:sin_TI*min}, $I^*$ as described in \eqref{eq:sine_I*} and \eqref{eq:sine_T} will be the optimal control for any spiking time $T\geq T^{I^*}_{min}$ if the bound satisfies $M\geq\omega/z_d$. A shorter spiking time $T\in[T^{M}_{min}, T^{I^*}_{min})$ is feasible but can not be achieved by $I^*$ alone. Suppose that $T\in[T^{M}_{min}, T^{I^*}_{min})$, then there exist two angles $\theta_1=\sin^{-1}[-2M\omega/(z_dM^2+cz_d)]$ and $\theta_2=\pi-\theta_1$ where $I^*$ meets the bound $M$, illustrated in Fig. 13(b). When $\theta\in(\theta_1,\theta_2)$, $I^*>M$ and we take $I(\theta)=M$ for $\theta\in[\theta_1,\theta_2]$. Then, from \eqref{eq:lambdac}, the Lagrange multiplier is $\lambda=(H-M^2)/(\omega+z_dM\sin\theta)$. This multiplier satisfies the adjoint equation \eqref{adj1}, therefore $I(\theta)=M$ is optimal for $\theta\in[\theta_1, \theta_2]$. Similarly, by symmetry, $I^*<-M$ when $\t\in[\t_3,\t_4]$, where $\theta_3=\pi+\theta_1$ and $\theta_4=2\pi-\theta_1$, if the desired spiking time is $T\in[T^{M}_{min}, T^{I^*}_{min})$. It can be easily shown by the same fashion that $I(\t)=-M$ is optimal in the interval $\t\in[\t_3,\t_4]$. Therefore, the minimum-power optimal control that spikes the neuron at $T\in[T^{M}_{min}, T^{I^*}_{min})$ can be characterized by four switchings between $I^*$ and $M$ as shown in \eqref{eq:I1*}.

\section{Morris-Lecar Model}
\label{appd:ML}
The dynamics of the Morris-Lecar neuron are described by two coupled dynamical equations 
\begin{align*}
	\dot{V} &= \frac{1}{C}\left[(I^b+I)+g_{Ca}m_{\infty}(V_{Ca}-V)+g_kw(V_k-V)+g_L(V_L-V)\right]\\
	\dot{w} &= \phi (\omega_{\infty}-w)/\tau_{w}(V)\\
	m_{\infty} &= 0.5[1+\tanh((V-V_1)/V_2)]\\
	\omega_{\infty}&=0.5[1+\tanh((V-V_3)/V_4)]\\
	\tau_{\omega}&=1/\cosh[(V-V_3)/(2V_4)].
\end{align*}
In Section \ref{sec:ML}, we consider the following parameter values
\begin{equation*}
    \begin{array}{lll}\phi=0.5, & I^b=0.09\ \mu A/cm^2, & V_1=-0.01\ mV \\
	V_2=0.15\ mV, & V_3=0.1\ mV, & V_4=0.145\ mV, \\
    g_{Ca}=1\  mS/cm^2, & V_k=-0.7\ mV, & V_L=-0.5\ mV,\\
    g_k=2\ mS/cm^2, & g_L=0.5\ mS/cm^2, & C=1\ \mu F/cm^2 \\
    Vca=1\ mV
	\end{array}
\end{equation*}

\section{}
\label{appd:HH}
The dynamics of the Hodgkin-Huxley neuron are described by a set of differential equations
\begin{align*}
	C\dot{V}-I &= -g_{Na}h(V-V_{Na})m^3-g_k(V-V_k)n^4-g_L(V-V_L)\\
	\dot{m} &= a_m(V)(1-m)-b_m(V)m\\
    \dot{h} &= a_h(V)(1-h)-b_h(V)h\\
    \dot{n} &= a_n(V)(1-n)-b_n(V)n\\
	a_m(V)&=0.1(V+40)/[1-\exp(-(V+40)/10)]\\
	b_m(V)&=4\exp[-(V+65)/18]\\
	a_h(V)&=0.07\exp[-(V+65)/20]\\
	b_h(V)&=1/(1+\exp[-(V+35)/10)]\\
	a_n(V)&=0.01(V+55)/[1-\exp(-(V+55)/10)]\\
	b_n(V)&=0.125\exp[-(V+65)/80].
 \end{align*}
In Section \ref{sec:HH}, we consider the following parameter values
$$\begin{array}{lll}V_{Na}=50\ mV, & V_k=-77\ mV, & v_L=-54.4\ mV, \\
	g_{Na}=120\  mS/cm^2,& g_k=36\ mS/cm^2, & g_L=0.3\ mS/cm^2,\\
	C=1\ \mu F/cm^2,& I= 10\ \mu A/cm^2. \end{array}$$

\section{Legendre pseudospectral method for optimal control of phase-reduced oscillators}
\label{sec:PSmethod}
The pseudospectral method is a spectral collocation method that was originally developed to solve partial differential equations, and has recently been adapted to solve optimal control problems \cite{Ross03,Li09, Li_JCP11, Li_PNAS11,Gong06,Elnagar95}. In this approach the differential equations that relate the states and the controls are discretized at specific collocation nodes, which results in a discrete optimization problem. All continuous-time functions are rescaled to the time domain of [-1,1] and expanded by an orthogonal polynomial basis based on a set of selected quadrature nodes \cite{Elnagar95}. Here, we use the Legendre-Gauss-Lobatto(LGL) nodes, and can then write the $N$th order interpolating approximations of the state and control functions
	\begin{align*}
		\theta(t)\approx I_N\theta(t)=&\sum_{k=0}^N\bar{\theta}_k\ell_k(t),\\
		I(t)\approx I_NI(t)=&\sum_{k=0}^N\bar{I}_k\ell_k(t),
	\end{align*}
	where
	$$\ell_k(t)=\prod_{i=0,i\neq k}^N\frac{t-t_i}{t_k-t_i},\quad k=0,1,\ldots,N,$$
	are the Lagrange polynomials with $\ell_k(t_i)=\delta_{ki}$, the Kronecker delta function. The derivative of $I_N\theta(t)$ at the LGL node $t_j$, $j=0,1,2,\ldots,N$ is then given by
	\begin{equation*}
	\frac{d}{dt}I_N\theta(t_j)=\sum_{k=0}^N\bar{\theta}_k\dot{\ell}_k(t_j)=\sum_{k=0}^ND_{jk}\bar{\theta}_k,
	\end{equation*}
	where $D_{jk}$ are the $jk^{th}$elements of the constant $(N+1)\times(N+1)$ differentiation matrix defined by
	\begin{equation*}
		D_{jk}=\left\{\begin{array}{ll}\frac{L_N(t_j)}{L_N(t_k)}\frac{1}{t_j-t_k} & j\neq k \\ -\frac{N(N+1)}{4} & j=k=0\\ \frac{N(N+1)}{4} & j=k=N \\ 0 & \mathrm{otherwise}.\end{array}\right.
	\end{equation*}
	The integral cost functional of the optimal control problem as in \eqref{eq:oc1} can be accurately approximated by the Gauss-Lobatto integration rule. Thus, the pseudospectral discretization of the optimal control problem \eqref{eq:oc1} gives rise to a nonlinear program of the form
	 \begin{align*}
	  \min_{\bar{I}_0\ldots \bar{I}_N} \quad & \frac{T}{2}\sum_{i=0}^N {\bar{I}}^2_iw_i,\\
	     {\rm s.t. } \quad & \sum_{k=0}^ND_{ik}\bar{\theta}_k=\frac{T}{2}\left[f(\bar{\theta}_i)+\bar{I}_ig(\bar{\theta}_i)\right],\\
	     &\sum_{k=0}^ND_{ik}\bar{p}_k =\frac{T}{2}\bar{I}_i, \\
	     &\bar{\theta}_0=0,  \quad \bar{\theta}_N=2\pi, \\
	     &\bar{p}_0=0,  \quad \bar{p}_N=0,\\
	     &|\bar{I}_i|\leq M,
	 \end{align*}
	 where $w_i$ are the LGL weights given by $ w_i=\frac{2}{N(N+1)}\frac{1}{(L_N(t_i))^2}$. %
	Solvers for this type of minimization problems are readily available and straightforward to implement. We approximate the problem using 151 nodes ($N=150$) and implement it in the AMPL language \cite{Gay93}. We use a third party nonlinear programming solver KNITRO \cite{Byrd99} to solve this optimization. This Legendre pseudospectral method provides a direct method to verify the analytical results presented in Section \ref{sec:mim_power_control}.

%
%

\ifCLASSOPTIONcaptionsoff
  \newpage
\fi




\bibliographystyle{IEEEtran}
\bibliography{IEEEabrv,charge_balanced}

\end{document}